\documentclass{gtart_h}  

\def\ifplaintex{\expandafter\ifx\csname documentclass\endcsname\relax}

\ifplaintex 
\else
\fi

\expandafter\ifx\csname beginpicture\endcsname\relax
\expandafter\ifx\csname documentclass\endcsname\relax
\input pictex \else
\input prepictex \input pictex \input postpictex \fi\fi

\def\gt{{\mathsurround=0pt\it $\cal G\mskip-2mu$eometry \&\ 
$\cal T\!\!$opology}}        

\def\gtp{{\mathsurround=0pt\it $\cal G\mskip-2mu$eometry \&\ 
$\cal T\!\!$opology $\cal P\!$ublications}}  

\def\lognumber#1{\def\thelognumber{#1}}
\def\volumenumber#1{\def\thevolumenumber{#1}}
\def\papernumber#1{\def\thepapernumber{#1}}
\def\volumeyear#1{\def\thevolumeyear{#1}}

\def\pagenumbers#1#2{\def\startpage{#1}\def\finishpage{#2}}
\def\published#1{\def\publishdate{#1}}
\def\proposed#1{\def\theproposer{#1}}
\def\seconded#1{\def\theseconders{#1}}
\def\received#1{\def\receiveddate{#1}}
\def\revised#1{\def\reviseddate{#1}}
\def\accepted#1{\def\accepteddate{#1}}
\def\asciititle#1{\def\theasciititle{#1}}

\def\asciiaddress#1{\def\theasciiaddress{#1}}
\def\asciiemail#1{\def\theasciiemail{#1}}

\long\def\asciiabstract#1{\long\def\theasciiabstract{#1}}
\def\asciikeywords#1{\def\theasciikeywords{#1}}

\let\\\par\let\thelognumber\relax
\let\thevolumenumber\relax\let\thepapernumber\relax
\let\thevolumeyear\relax\let\thesamplenumber\relax\let\startpage\relax
\let\finishpage\relax\let\publishdate\relax\let\receiveddate\relax
\let\reviseddate\relax\let\accepteddate\relax\let\theasciititle\relax
\let\theasciiauthors\relax\let\theasciiaddress\relax
\let\theasciiabstract\relax\let\theasciikeywords\relax
\let\theasciiemail\relax\let\theshortauthors\relax\let\theshorttitle\relax

\long\def\maketitlep{   

\count0=\startpage

\gt\hfill      
\beginpicture
\setcoordinatesystem units <0.33truein, 0.33truein> point at 2.2 0.9
\setplotsymbol ({$\cal G$})
\plotsymbolspacing=9truept
\circulararc 315 degrees from 0 1 center at 0 0
\setplotsymbol ({$\cal T$})
\circulararc 315 degrees from 1 -1 center at 1 0
\endpicture
\break
{\small\ifx\thesamplenumber\relax 
Volume \else Sample
\fi\thevolumenumber\ (\thevolumeyear)
\startpage--\finishpage\nl
Published: \publishdate}
\vglue 0.5truein plus 0.4fil minus 0.1truein

{\parskip=0pt\leftskip 0pt plus 1fil\def\\{\par\smallskip}{\ifplaintex\large
\else\Large\fi\bf\thetitle}\par\medskip}   

\vglue 0pt plus 0.1fil 

{\parskip=0pt\leftskip 0pt plus 1fil\def\\{\par}{\sc\theauthors}
\par\medskip}

\vglue 0pt plus 0.1fil 

{\small\parskip=0pt\let\newline\\
{\leftskip 0pt plus 1fil\def\\{\par}{\sl\theaddress}\par}
\expandafter\ifx\theemail\relax    
\relax\else\vglue 5pt plus 0.02fil minus 2pt\def\\{\stdspace{\rm 
and}\stdspace} 
\cl{Email:\stdspace\tt\theemail}\fi
\ifx\theurl\relax                  
\relax\else\vglue 5pt plus 0.02fil minus 2pt\def\\{\stdspace{\rm 
and}\stdspace}
\cl{URL:\stdspace\tt\theurl}\fi\par}

\vglue 7pt plus 0.3fil minus 3pt

{\bf Abstract}
\vglue 5pt plus 0.1fil minus 2pt

\theabstract

\vglue 7pt plus 0.3fil minus 3pt

{\bf AMS Classification numbers}\quad Primary:\quad \theprimaryclass

Secondary:\quad \thesecondaryclass

\vglue 5pt plus 0.3fil minus 2pt

{\bf Keywords:}\quad \thekeywords

\vglue 10pt plus 0.5fil minus 5pt

{\small  Proposed: \theproposer\hfill Received: \receiveddate\nl
Seconded: \theseconders\hfill 
\ifx\reviseddate\relax                         
Accepted: \accepteddate                        
\else
Revised: \reviseddate                          
\fi}
\eject
}       

\let\maketitlepage\maketitlep
\let\maketitle\maketitlepage

\font\phead=cmsl9 scaled 950
\font=cmsl9 scaled 1050
\font\pnum=cmbx10 scaled 913
\font\lnum=cmbx10 
\font\pfoot=cmsl9 scaled 950
\font=cmsl9 scaled 1050
\ifplaintex
\headline{\vbox to 0pt{\vskip -4.5mm\line{\small\phead\ifnum
\count0=\startpage ISSN 1364-0380 (on line)
1465-3060 (printed) \hfill {\pnum\folio}\else\ifodd\count0\def\\{ }%
\ifx\theshorttitle\relax\thetitle\else\theshorttitle\fi\hfill{\pnum\folio}
\else\def\\{ and }{\pnum\folio}\hfill\ifx\theshortauthors\relax\theauthors
\else\theshortauthors\fi\fi\fi}\vss}}
\footline{\vbox to 0pt{\vglue 0mm\line{\small\pfoot\ifnum\count0=\startpage
\copyright\ \gtp\hfill\else
\gt, Volume \thevolumenumber\ (\thevolumeyear)\hfill\fi}\vss
}}
\else
\makeatletter
\def\@oddhead{{\small\lhead\ifnum\count0=\startpage ISSN 1364-0380 (on line)
1465-3060 (printed) \hfill {\lnum\number\count0}\else\ifodd\count0
\def\\{ }\ifx\theshorttitle\relax \thetitle \else\theshorttitle\fi\hfill
{\lnum\number\count0}\else\def\\{ and }{\lnum\number\count0}
\hfill\ifx\theshortauthors\relax 
\theauthors\else\theshortauthors\fi\fi\fi}}\def\@evenhead{\@oddhead}
\def\@oddfoot{\small\lfoot\ifnum\count0=\startpage\copyright\ \gtp\hfill\else
\gt, Volume \thevolumenumber\ (\thevolumeyear)\hfill\fi}
\def\@evenfoot{\@oddfoot}
\makeatother
\fi

\newwrite\gtoutfile
\long\gdef\makeheadfile{  
{\def\\{, }\def\s{ }
\immediate\openout\gtoutfile head.xxx
\immediate\write\gtoutfile{Proxy-for: \ifx\theasciiauthors\relax
\theauthors\else\theasciiauthors\fi\s<\ifx\theasciiemail\relax\theemail\else\theasciiemail\fi>}
\immediate\write\gtoutfile{\noexpand\\}
\immediate\write\gtoutfile{Authors: \ifx\theasciiauthors\relax
\theauthors\else\theasciiauthors\fi}
{\def\\{ }\immediate\write\gtoutfile{Title: \ifx\theasciititle\relax
\thetitle\else\theasciititle\fi}}
\immediate\write\gtoutfile{Subj-class: GT or SG or MG etc}
\immediate\write\gtoutfile{MSC-class: \theprimaryclass\ifx\thesecondaryclass\relax\else, \thesecondaryclass\fi}
\immediate\write\gtoutfile{Journal-ref: Geom. Topol. \thevolumenumber
(\thevolumeyear) \startpage-\finishpage}
\immediate\write\gtoutfile{Comments: Published by Geometry and Topology at}
\immediate\write\gtoutfile{\s\s http://www.maths.warwick.ac.uk/gt/GTVol\thevolumenumber/paper\thepapernumber.abs.html}
\immediate\write\gtoutfile{\noexpand\\}
\immediate\write\gtoutfile{}
\ifx\theasciiabstract\relax
\immediate\write\gtoutfile{\theabstract}\else
\immediate\write\gtoutfile{\theasciiabstract}\fi
\immediate\write\gtoutfile{}
\immediate\write\gtoutfile{\noexpand\\}
\immediate\write\gtoutfile{}
\immediate\closeout\gtoutfile}}  

\def\maketitlepage{\maketitlep\makeheadfile}
\let\maketitle\maketitlepage

\lognumber{452}
\received{06 May 2004}
\volumenumber{9}\papernumber{45}\volumeyear{2005}
\pagenumbers{1977}{2012}   
\revised{16 October 2005}
\published{17 October 2005}
\accepted{24 April 2005}
\proposed{Ronald Stern}
\seconded{Ronald Fintushel, Robion Kirby}

\usepackage{amssymb,amsmath}

\allowdisplaybreaks

\newtheorem{theorem}{Theorem}[section]
\newtheorem{prop}[theorem]{Proposition}
\newtheorem{lemma}[theorem]{Lemma}

\theoremstyle{remark}

\newtheorem{rem}[theorem]{Remark}

\def\Z{{\mathbb Z}}

\def\M{{\cal M}}
\def\CM{\overline{\cal M}}

\def\H{{\cal H}}

\def\a{\alpha}
\def\b{\beta}
\def\g{\gamma}

\def\s{{\rm s}}
\def\e0{\epsilon_{0}}

\begin{document}

\title{Yau--Zaslow formula on K3 surfaces\\for non-primitive classes}
\asciititle{Yau-Zaslow formula on K3 surfaces for non-primitive classes}
\author{Junho Lee\\Naichung Conan Leung}
\address{Department of Mathematical Sciences, Seoul National 
University San56-1\\Shinrim-dong Kwanak-gu, Seoul 151-747, Korea }
\secondaddress{Institute of Mathematical Sciences\\The 
Chinese University of Hong Kong, Shatin, NT, Hong Kong}
\asciiaddress{Department of Mathematical Sciences, Seoul National 
University San56-1\\Shinrim-dong Kwanak-gu, Seoul 151-747, 
Korea\\and\\Institute of Mathematical Sciences\\The 
Chinese University of Hong Kong, Shatin, NT, Hong Kong}
\asciiemail{leejunho@msu.edu, leung@ims.cuhk.edu.hk}
\gtemail{\mailto{leejunho@msu.edu}, \mailto{leung@ims.cuhk.edu.hk}}

\begin{abstract}   
We compute the genus zero family Gromov--Witten invariants for K3
surfaces using the topological recursion formula and the
symplectic sum formula for a degeneration of elliptic K3 surfaces.
In particular we verify the Yau--Zaslow formula for non-primitive
classes of index two.
\end{abstract}
\asciiabstract{%
We compute the genus zero family Gromov-Witten invariants for K3
surfaces using the topological recursion formula and the symplectic
sum formula for a degeneration of elliptic K3 surfaces.  In particular
we verify the Yau-Zaslow formula for non-primitive classes of index
two.}

\asciikeywords{Family Gromov-Witten invariants, Yau-Zaslow formula,
          symplectic sum formula, topological recursion relation, K3 surface }


\maketitlepage    


\setcounter{equation}{0}
\section{Introduction}

Let $N\left( d,r\right) $ be the number of rational curves in K3
surfaces $X$ that represent a homology class $A\in H_{2}\left(
X,{\Bbb Z}\right) $ of self-intersection $A^{2}=2d-2$ and of
index\footnote{The index of $A$ is the largest positive integer
$r$ such that $r^{-1}A$ is integral. An index one class is called
primitive.} $r$. Yau and Zaslow \cite{yz} give an ingenious
heuristic argument to compute the generating function for
primitive classes and they also expected that the same formula holds
true for classes of arbitrary index. More precisely the Yau--Zaslow
conjectural formula says that, for any positive integer $r$, we
have
\begin{equation}\label{0.1}
 \sum_{d\geq 0}N\left( d,r\right) \,t^{d}\,=\,
 \prod_{d\geq 1}\Big(\frac{1}{1-t^{d}}\Big)^{24}.
\end{equation}
Their original approach was pursued by Beauville \cite{b}, Chen
\cite{c} and Li \cite{li}. In \cite{bl,bl2}, Bryan and the second
author showed that $N\left( d,r\right) $ can be computed in terms
of the twistor family Gromov--Witten invariants of the K3 surfaces.
They also proved the Yau--Zaslow formula for primitive classes in
K3 surfaces and its higher genera generalization.

In \cite{l2}, the first author reproved the Yau--Zaslow formula (\ref{0.1})
for primitive classes and its higher genera generalization
using $p_{g}$--dimensional
family Gromov--Witten invariants defined in \cite{l1} --- following the approach
of \cite{ip3}, he computed those invariants by relating the TRR
(topological recursion relation) and the symplectic sum formula of
\cite{ip3} for a suitable degeneration of an elliptic K3 surface.
In this article we explain how to use the same approach to compute
the $p_{g}$--dimensional family Gromov--Witten invariants for
non-primitive classes. In particular {\em we verify the Yau--Zaslow
formula for non-primitive classes of index two}. At present, it is
not easy to use this approach to handle classes of higher indexes,
for example, we do not know how to handle relative invariants with
multiplicity greater than 2. An analogous problem for the
Seiberg--Witten invariants was studied by Liu \cite{l}.

Notice that $N\left( d,2\right) $ is different from the family
Gromov--Witten invariant $GW^{\H}_{A,0}$ due to the multiple cover
contributions, as it was explained by Gathmann in \cite{ga}. This
is because the family Gromov--Witten invariants count holomorphic
maps and each rational curve $C$ representing the primitive class
$A/2$ contributes $1/2^{3}$ to $GW^{\H}_{A,0}$, however, the
multiple curve $2C$ contributes zero to $N\left( d,2\right) $
because it has negative genus. As a result, the Yau--Zaslow formula
for non-primitive classes of index two follows directly from the
following theorem.

\begin{theorem}
\label{T:Main} Let $X$ be a K3 surface and $A/2\in
H_2\big(X;{\mathbb Z}\big)$ be a primitive class. Then, the genus
$g=0$ family GW invariant of $X$ for the class $A$ is given by
\begin{equation*}
  GW_{A,0}^{{\cal H}}\,-\,GW_{B,0}^{{\cal H}}\,=\,
 \Big(\frac{1}{2}\Big)^{3}\,GW_{A/2,0}^{{\cal H}}  
\end{equation*}
where $B$ is any primitive class with $B^2=A^2$.
\end{theorem}

In the sequel \cite{lle}, we apply the same technique to enumerate
the number of elliptic curves representing non-primitive classes
of index two in K3 surfaces.

The construction of family GW invariants is briefly
described in Section 2. We outline the proof of
Theorem~\ref{T:Main} in Section 3. This proof follows the elegant
argument used by Ionel and Parker to compute the GW invariants of
$E(0)$ \cite{ip3}. It involves computing the generating functions
for the invariants in two different ways, first using the TRR
formula, and second using the symplectic sum formula. Section 4
gives the sum formulas of the symplectic sum of $E(2)$ with $E(0)$
along a fixed fiber. The sum formulas yield relations of family
invariants of $E(2)$ and relative invariants of $E(0)$. We compute
those relative invariants of $E(0)$ in Section 5--7.

\medskip
{\bf Acknowledgments}\qua The first author would like to thank Thomas
Parker for his extremely helpful discussions and he is also grateful
to Eleny Ionel, Bumsig Kim, Tian Jun Li, Ionut Ciocan-Fontaine, Jun Li
and Dihua Jiang for their useful comments. In addition, the first
author wish to thank Ronald Fintushel for his interest in this work
and especially for his encouragement.  We are also very grateful for
the referee's useful comments. The second author is partially
supported by NSF/DMS-0103355, CUHK/2060275, and CUHK/2160256.


\setcounter{equation}{0}
\section{Family GW invariants of K3 surfaces}
\label{section1}

This section briefly describes family GW invariants of K3
surfaces. We first give the definition of family GW invariants of
K\"{a}hler surfaces with $p_{g}=\mbox{dim}H^{2,0}\geq 1$ defined
in \cite{l1}. Fix a compact K\"{a}hler surface $(X,J)$ and choose
the $2p_{g}$--dimensional parameter space
\begin{equation*}
 \H\,=\,\mbox{Re}\,\big(H^{0,2}(X)\oplus H^{0,2}(X)\big).
\end{equation*}
Using the K\"{a}hler metric, each $\a\in\H$ defines an
endomorphism $K_{\a}$ of $TX$ by the equation
\begin{equation*}
 \langle u,K_{\a}v \rangle\,=\,\a(u,v).
\end{equation*}
Since $Id+JK_{\a}$ is invertible for each $\a\in \H$, the equation
\begin{equation*}
 J_{\a}\,=\,(Id\,+\,J\,K_{\a})^{-1}J\,(Id\,+\,J\,K_{\a})
\end{equation*}
defines a family of almost complex structure on $X$ parameterized
by $\a$ in the $2p_{g}$--dimensional linear space $\H$. The family
GW invariants are defined, in the same manner as the ordinary GW
invariants \cite{rt,lt}, but using the moduli space of stable
$(J,\alpha)$--holomorphic maps $(f,\alpha)$:
\begin{equation}
 \CM^{\H}_{g,k}(X,A,J) = \{\ (f,\alpha)\ |\ \overline{
 \partial}_{J_{\alpha}} = 0\,,\,\, \alpha\in{\cal H}\,,\, [f]=A\in H_{2}(X;
 {\mathbb Z})\, \}.  \label{modulispace}
\end{equation}
For each stable $(J,\a)$ holomorphic map $ f\co (C,j)\rightarrow X$
of genus $g$ with $k$--marked points, collapsing unstable
components of the domain determines a point in the Deligne--Mumford
space $\CM_{g,k}$ and evaluation of marked points determines a
point in $X^{k}$. Thus we have  a map
\begin{equation}\label{st-ev}
  st\times ev\ \co \ \CM^{\H}_{g,k}(X,A,J)\ \to\
  \CM_{g,k}\times X^{k}
\end{equation}
where $st$ and $ev$ denote the stabilization map and the
evaluation map, respectively. If  the space
(\ref{modulispace}) is compact,  it  carries a fundamental
homology class
$$[\,\CM^{\H}_{g,k}(X,A,J)\,]$$
which we can push
forward by the map (\ref{st-ev}) to obtain a  homology
class
\begin{equation*}
 (st\times ev)_{*}\,\big[\,\CM^{\H}_{g,k}(X,A,J)\,]\, \in\, H_{2r}\big(\CM_{g,k}\times X^{k};{\Bbb Q}\big)
\end{equation*}
where $r=-K\cdot A+(g-1)+k+p_{g}$ and $K$ is the canonical class
of $X$. Then, the family GW invariants are defined by
\begin{align*}
 &GW^{\H}_{g,k}(X,A,J)\big(\kappa;\b_{1},\cdots,\b_{k}\big) \\
 &=\, (st\times ev)_{*}\,
 \big[\,\CM^{\H}_{g,k}(X,A,J)\,]\,\cap\,
 \big( \kappa^{*}\,\cup\,\b_{1}^{*}\,\cup\,\cdots\,\cup\,\b_{k}^{*}\big)
\end{align*}
where $\kappa^{*}$ and $\b_{i}^{*}$ are Poincar\'{e} dual of
$\kappa\in H_{*} \big(\CM_{g,k};{\Bbb Q}\big)$ and $\b_{i}\in
H_{*}\big(X^{k};{\Bbb Q}\big)$, respectively.

When $X$ is a K3 surface, the family GW invariants reduce to the
invariants defined by Bryan and Leung \cite{bl} using the twistor
family. In particular, (i) they are independent of complex
structures and (ii) for any two homology classes $A$ and $B$ of
the same index with $A^{2}=B^{2}$, there is an orientation
preserving diffeomorphism $h\co X\to X$ such that $h_{*}A=B$ and
\begin{equation}  \label{T:BL}
 GW_{g,k}^{\H}(X,A)(\kappa;\b_{1},\cdots,\b_{k})\,=\,
 GW_{g,k}^{\H}(X,B)(\kappa;h_{*}\b_{1},\cdots,h_{*}\b_{k}).
\end{equation}
Below, we will often write the family GW invariants of $K3$
surfaces as simply
\begin{equation*}
 GW^{\H}_{A,g}(X)\big(\kappa;\b_{1},\cdots,\b_{k}\big)
 \ \ \ \ \ \mbox{or}\ \ \ \ \
 GW^{\H}_{A,g}\big(\kappa;\b_{1},\cdots,\b_{k}\big).
\end{equation*}
By dimension count, this invariant vanishes unless
\begin{equation*}
 \mbox{deg}(\kappa^{*})\,+\,\sum\mbox{deg}(\b_{i}^{*})
 \,=\,2(g+k).
\end{equation*}

Let $E(2)\to {\Bbb P}^{1}$ be an elliptic K3 surface with a
section of self intersection number $-2$. Denote by $S$ and $F$
the section class and the fiber class, respectively. It then
follows from (\ref{T:BL}) that for any class $A$ of index 2 with
$(A/2)^{2}=2d-2$, we have
\begin{equation*}
GW^{\H}_{A,0}\,=\,GW^{\H}_{2(S+dF),0}.
\end{equation*}
These are the invariants we aim to compute. We will compute them
following a similar approach of \cite{ip3,l2} ---  relating the
genus 1 TRR Formula and the Symplectic Sum Formula \cite{ip3}.

The Sum Formula yields relations between GW invariants and
relative GW invariants of \cite{ip2}. One can derive a family
version of  the Sum Formula for the cases of elliptic
surfaces (cf \cite{l2}). Here, we introduce relative family
invariants of $E(2)$  and describe the extension of the Sum
Formula in Section 3.

First, we define relative family invariants of $E(2)$ for the
classes $2S+dF$, $d\in {\Bbb Z}$. Let $V\cong T^{2}$ be a smooth
fiber of $E(2)\to {\Bbb P}^{1}$ and choose a smooth bump function
$\mu$ that vanishes in a small $\delta$--neighborhood of $V$ and is
1 everywhere outside of a $2\delta$--neighborhood of $V$. Replacing
$\a\in\H$ by $\mu\a$ and following the construction of relative
invariants in \cite{ip2}, one can define the moduli space of
`$V$--regular' $(J,\mu\a)$--holomorphic maps $(f,\mu\a)$
\begin{equation}  \label{rmoduli}
 \M^{\H,V}_{g,k,s}\big(2S+dF\big)
\end{equation}
where $s=(s_{1},\cdots,s_{l})$ is a multiplicity vector and
$f^{-1}(V)$ consists of marked points $p_{j}$, $k+1\leq j\leq
k+l$, each with the contact order of $f$ with $V$ at $p_{j}$ being
$s_{j}$. Since each $s_{j}\geq 1$ and $(2S+dF)\cdot [V]=2$, the
multiplicity vector $s$ is either $(1,1)$ or $ (2)$. This moduli
space also comes with a map
\begin{equation}\label{st-ev-h}
  st\times ev\times h\,:\
  \M^{\H,V}_{g,k,s}(S+dF)\
  \to\
  \CM_{g,k+l}\times E(2)^{k}  \times V^{l}
\end{equation}
where  $ev$ is the evaluation map of first $k$ marked points into
$E(2)^{k}$ and $h$ is the evaluation map of last $l$ marked points
into $V^{l}$. The moduli space (\ref{rmoduli}) is compact (cf
Section 6 of \cite{l2}) and hence carries a fundamental class
$$ [\M^{V}_{g,k,s}(S+dF)] $$
which we can push forward by the map
(\ref{st-ev-h}) to obtain a homology class
\begin{equation}  \label{rel-fc}
(st\times ev\times h)_{*}\,
 \big[\M^{\H,V}_{g,k,s}\big(2S+dF\big)\big] \in H_{2r}\big(\overline{{\cal M}}
_{g,k+l}\times E(2)^{k}\times V^{l};{\Bbb Q}\big)
\end{equation}
where $r=g+k+l-2$. The relative family invariants of $(E(2),V)$ is
then defined as
\begin{align}\label{RI-V}
 &GW^{V}_{2S+dF,g,s}
 \big(\kappa;\b_{1},\cdots,\b_{k};C_{\g_{1}\cdots\g_{l}}\big)
 \notag \\ &=\
 (st\times ev\times h)_{*}\,
 \big[\M^{\H,V}_{g,k,s}\big(2S+dF\big)\big]\,\cap\,
 \big(\kappa^{*} \cup \b^{*}\cup \g^{*}\big)
\end{align}
where $\b^{*}=\b_{1}^{*}\cup\cdots\cup\b_{k}^{*}$,
$\g^{*}=\g_{1}^{*}\cup\cdots\cup\g_{l}^{*}$, and $\g_{j}^{*}$ is
the Poincar\'{e} dual of $\g_{j}\in H_{*}(V;{\Bbb Z})$ in $V$.

Similarly, we can define relative family invariants of $E(2)$ for
the classes $S+(2d-3)F$, $d\in {\Bbb Z}$. In this case, we choose
a symplectic submanifold $U\cong T^{2}$ of $E(2)$ that represents
the class $2F$. Repeating the same arguments as above then gives
relative family invariants of $(E(2),U)$
\begin{equation}\label{RI-U}
 GW_{S+(2d-3)F,g,s}^{U}(\kappa;\b_1,\cdots,\b_k;C_{\g_{1}\cdots\g_{l}})
\end{equation}
where $s$ is also either (1,1) or (2) since each $s_{j}\geq 1$ and
$(S+(2d-3)F)\cdot [U]=2$.

The invariant (\ref{RI-V}) (resp.\ (\ref{RI-U})) counts the
oriented number of genus
 $g$  $V$--regular (resp. $U$--regular)  $(J,\mu\a)$--holomorphic maps $(f,\mu\a)\co C\to E(2)$,
 representing the homology class $2S+dF$ (resp.\ $S+(2d-3)F$),
 with $C\in K$ and $f(x_{i})\in A_{i}$ such that these have a contact of order
 $s_{j}$ with $V$ (resp.\ $U$) along fixed representatives $G_{j}$ of $\g_{j}$ in $V$ (resp.\ $U$) where
 $K$ and $A_{i}$ are representatives of $\kappa$ and $\b_{i}$.
In particular, both relative invariants (\ref{RI-V}) and
(\ref{RI-U}) have the same (formal) dimension and thus vanish
unless
\begin{equation*}
 \mbox{deg}(\kappa^{*})\,+\,\sum\mbox{deg}(\b_{i}^{*})
 \,+\,\sum\mbox{deg}(\g_{j}^{*})\,=\,2\big(g+k+l-2\big).
\end{equation*}


\setcounter{equation}{0}
\section{Outline of computations}
\label{section2}

Our goal is to compute the $g=0$ family invariants of $K3$
surfaces for the classes $A$ of index 2. By (\ref{T:BL}), it
suffices to compute the invariants of $E(2)$ for the classes
$2(S+dF)$. For convenience we assemble them in the generating
functions
\begin{equation}\label{def-M}
 M_{g}(\,\cdot\,)(t) \,=\,
 \sum\, GW^{{\cal H}}_{2S+dF,g}\big(\,\cdot\,\big) \,t^{d}\,.
\end{equation}
We further introduce generating functions for invariants of
primitive classes, by the formula
\begin{align}
 N_{g}(\,\cdot\,)(t)\,&=\,\sum\,
 GW^{\H}_{S+dF,g}\big(\,\cdot\,\big)\,t^{d}\,,\notag \\
  \vspace{3cm}
 P_{g}(\,\cdot\,)(t)\,&=\,\sum\,
  GW^{\H}_{S+(2d-3)F,g}\big(\,\cdot\,\big)\,t^{d}.
  \label{def-P}
\end{align}
It then follows from  (\ref{T:BL}) that
\begin{equation*}
 M_{0}(t)\,-\,P_{0}(t)\,=\,
 \sum_{d\geq 0}\big(\,
 GW^{\H}_{2(S+dF),0}\,-\,GW^{\H}_{S+(4d-3)F,0}\,\big)\,t^{2d}
\end{equation*}
since both $2S+dF$ and $S+(2d-3)F$ are primitive with the same
square when $d$ is odd. In this and the following four sections we
will show:

\begin{prop}
\label{Main-Prop}\ \ \ 
 $M_{0}(t) - P_{0}(t)\,=\,\Big(\dfrac{1}{2}\Big)^{3}\,N_{0}(t^{2}).$
\end{prop}

Let $A$ be any class of index 2 and $B$ be any primitive classes
such that $A^{2}=B^{2}=4(2d-2)$. Proposition~\ref{Main-Prop} and
(\ref{T:BL}) then imply
\begin{align*}
 GW_{A,0}^{{\cal H}}\,-\,GW_{B,0}^{{\cal H}}\ &=\
 GW^{\H}_{2(S+dF),0}-GW^{\H}_{S+(4d-3)F,0}\\
 &=\ \big(\tfrac{1}{2}\big)^{3}\,GW^{\H}_{S+dF,0}\ =\
 \big(\tfrac{1}{2}\big)^{3}\,GW^{\H}_{A/2,0}
\end{align*}
and hence prove Theorem~\ref{T:Main} of the introduction.

Below, we outline 	how we prove Proposition~\ref{Main-Prop}.

Let $\varphi_{i}$ be  the first Chern class of the line bundle
${\cal L}_{i}^{\H}\to \CM_{g,k}^{\H}\big(X,A\big)$ whose geometric
fiber at the point $\big(C;x_{1},\cdots,x_{k},f,\alpha\big)$ is
$T^{*}_{x_{i}}C$. Similarly as for the ordinary GW invariants, one
can use $\varphi_{i}$ to impose descendent constraints on family
invariants as follows:
\begin{align*}
 &GW_{g,k}^{\H}(X,A)(\tau_{m_{1}}(\b_{1}),\cdots,\tau_{m_{k}}(\b_{k}))\\
 &=\
 (st\times ev)_{*}\big(\,\big[\,\CM^{\H}_{g,k}(X,A)\,]\,\cap\,
  \varphi_{1}^{m_{1}}\cdots\varphi_{k}^{m_{k}}\,\big)\ \cap\
 (\,\b_{1}^*\cup\cdots\cup\b_{k}^*\,).
\end{align*}
If the constraint $\tau_{m_{i}}(\b_{i})$ repeats $n$ times and
$\mbox{deg}(\b^*_{i})$ is even, we will use the notation
$\tau_{m_{i}}(\b_{i})^{n}$.

Recall that $V\cong T^{2}$ is a fixed
smooth fiber of $E(2)\to {\Bbb P}^{2}$. To save notation, we
denote by $F$ the fundamental class of $V$. Introduce a generating
functions for the relative invariants of $(E(2),V)$, by the
formula
\begin{equation}\label{def-M^V}
 M^{V}_{1,(2)}(t) \,=\,
 \sum\,GW_{2S+dF,1,(2)}^{V}\big(C_{F}\big)\,t^{d}\,.
\end{equation}
As in Proposition 3.1 of \cite{l2}, we can combine the $g=1$ TRR
formula with the composition law (Proposition 3.7 of \cite{l1})
to have
\begin{equation}
 M_{1}\big(\tau(F)\big)\,=\,
 \tfrac{1}{3}\,t\,M^{\prime}_{0} \,- \,\tfrac{2}{3}
 \,M_{0}.  \label{TRRg=1}
\end{equation}
 Then, in
Proposition~\ref{P:sum-formula} we apply the Symplectic Sum
Formula of  \cite{ip3}  to obtain
\begin{align}
 M_{1}\big(\tau(F)\big)\,&=\, M^{V}_{1,(2)} \,+\, 4\,G_{2}\,M_{0},
 \label{Sum-1}\\
 \vspace{3cm}
 M_{2}\big(\tau(F)^2\big)\,-2\,M_{1}\big(pt\big)\,& =\,
 20\,G_{2}\,M^{V}_{1,(2)}\, +\,\left(16\,G^{2}_{2}\,
 +\,8\,t\,G^{\prime}_{2}\right)M_{0}  \label{Sum-2}
\end{align}
where  $G_{2}(t)$ is the Eisenstein series of weight 2, namely
\begin{equation*}
 G_{2}(t)\,=\, \sum_{d\geq 0}\sigma(d)\,t^{d} \ \ \ \mbox{where}\ \
 \  \sigma(d)\,=\,\sum_{k|d}k,\ d\geq 1\ \ \ \mbox{and}\ \ \ \sigma(0)\,=\,- \tfrac{1}{24}.
\end{equation*}
Now, eliminate $M^{V}_{1,(2)}$ in (\ref{Sum-2}) by using
(\ref{TRRg=1}) and (\ref{Sum-1}) to obtain
\begin{equation}  \label{ODE-1}
M_{2}\big(\tau(F)^2\big)\,
 -2\,M_{1}\big(pt\big)\,=\, \tfrac{20}{3}
 \,G_{2}\,t\,M^{\prime}_{0}\,-\, \left(64\,G_{2}^{2}\,+\, \tfrac{40}{3}
 \,G_{2}\,-\,8\,t\,G_{2}^{\prime}\right)M_{0}\,.
\end{equation}

Recall that $U\cong T^{2}$ is a fixed symplectic submanifold of
$E(2)$ that represents the class $2F$.  Without any further
confusion, we will also denote by $F$ the fundamental class of
$U\cong T^2$. Introduce a generating function for the relative
invariants of $(E(2),U)$, by the formula
\begin{equation}\label{def-P^U}
 P^{U}_{1,(2)}(t)\,=\,
 \sum\,GW_{S+(2d-3)F,1,(2)}^{U}\big(C_{F}\big)\,t^{d}\,.
\end{equation}
The genus 1 TRR formula gives a formula like (\ref{TRRg=1})
\begin{equation}
 P_{1}\big(\tau(2F)\big)\,=\,
 \tfrac{1}{3}\,t\,P^{\prime}_{0} \,- \,\tfrac{2}{3}
 \,P_{0}\,.  \label{TRRg=1-P}
\end{equation}
Then, in Proposition~\ref{P:sum-formula-U} we apply the sum
formula to have formulas like (\ref{Sum-1}) and (\ref{Sum-2})
\begin{align}
 P_{1}\big(\tau(2F)\big)\,&=\, P^{U}_{1,(2)} \,+\, 4\,G_{2}\,P_{0},
 \label{Sum-1-P}\\
 \vspace{3cm}
 P_{2}\big(\tau(2F)^{2}\big)\,-2\,P_{1}\big(pt\big)\,& =\,
 20\,G_{2}\,P^{U}_{1,(2)}\, +\,\left(16\,G^{2}_{2}\,
 +\,8\,t\,G^{\prime}_{2}\right)P_{0}  \label{Sum-2-P}
\end{align}
Similarly, as above, equations (\ref{TRRg=1-P}),
(\ref{Sum-1-P}) and (\ref{Sum-2-P}) give
\begin{equation}  \label{ODE-1-P}
 P_{2}\big(\tau(2F)^{2}\big)\,
 -2\,P_{1}\big(pt\big)\,=\, \tfrac{20}{3}
 \,G_{2}\,t\,P^{\prime}_{0}\,-\, \left(64\,G_{2}^{2}\,+\, \tfrac{40}{3}
 \,G_{2}\,-\,8\,t\,G_{2}^{\prime}\right)P_{0}\,.
\end{equation}
Note that the equations (\ref{ODE-1}) and (\ref{ODE-1-P}) have the
same coefficients. This is true because all coefficients of TRR
and sum formula depends only on the topological quantities
\begin{equation*}
 (2S+dF)^{2}\,=\,(S+(2d-3)F)^{2},\ \ 
 (2S+dF)\cdot F\,=\,(S+(2d-3)F)\cdot 2F,\ \ 
 F^{2}\,=\,(2F)^{2}.
\end{equation*}
Hence, (\ref{ODE-1}) and (\ref{ODE-1-P}) give
\begin{align}
 &3\,\big[\,M_{2}\big(\tau(F)^{2}\big)\,-\,
        P_{2}\big(\tau(2F)^{2}\big)\,\big]
  \,-\,
 6\,\big[\,M_{1}\big(pt\big)\,-\,
           P_{1}\big(pt\big)\,]\notag \\ \vspace{1cm}
 &=\, 20
 \,G_{2}\,t\,\big(\,M_{0}\,-\,P_{0}\,\big)^{\prime}\,-\,
 \left(192\,G_{2}^{2}\,+\, 40
 \,G_{2}\,-\,24\,t\,G_{2}^{\prime}\right)
 \big(\,M_{0}\,-\,P_{0}\,\big)\,. \label{ODE-1-1}
\end{align}
Note that by (\ref{T:BL}) both generating functions
 $M_{1}\big(pt\big)-P_{1}\big(pt\big)$ and
 $M_{0}-P_{0}$ have no odd terms.
One can also show that the generating function
 $M_{2}\big(\tau(F)^{2}\big) -
        P_{2}\big(\tau(2F)^{2}\big)$ has no odd
terms (see \cite{lle}). Consequently, comparing odd terms of both
sides of (\ref{ODE-1-1}) gives the first order ODE
\begin{equation}  \label{ODE-3}
 0 =\, 20\,G_{o}\,t\,
 \big(\,M_{0}-P_{0}\,\big)^{\prime}\,-\,
 \left(384\,G_{e}\,G_{o}\,+\,40\,G_{o}\,
 -\,24\,t\,G_{o}^{\prime}\right)\big(\,M_{0}-P_{0}\,\big)\,
\end{equation}
where $G_{e}(t)$ (resp. $G_{o}(t)$) is the sum of all
even (resp. odd) terms of $G_{2}(t)$.

On the other hand, it follows from the equation (2.7) of \cite{l2}
that
\begin{equation}  \label{primitive}
 t\,\frac{d}{dt}N_{0}(t^{2})\,=\,
 2\,t^{2}\,N_{0}^{\prime}(t^{2})\,=\,
 48\,G_{2}(t^{2})\,N_{0}(t^{2})\,+\,2\,N_{0}(t^{2}).
\end{equation}
Combining (\ref{primitive}) with the following relation of certain
quasi-modular forms
\begin{equation}  \label{rel-qmod}
 4\,t^{2}\,G_{2}^{\prime}(t^{2})\,=\,32\,G_{2}^{2}(t^{2})\,-\,
 40\,G_{2}(t)\,G_{2}(t^{2})\,+\,8\,G_{2}^{2}(t)\,-\,t\,G_{2}^{\prime}(t)
\end{equation}
(see Lemma~\ref{L:qmf}) says that $N_{0}(t^{2})$ also satisfies
the same ODE  (\ref{ODE-3}). Since the initial conditions are
$N_{0}(0)=1$\ (cf \cite{l2}) and $M_{0}(0)-P_{0}(0)=(1/2)^{3}$\
(cf \cite{ga}), we can conclude that
\begin{equation*}
\label{mul-con}
M_{0}(t)-P_{0}(t)\,=\,\Big(\tfrac{1}{2}\Big)^{3}\,N_{0}(t^{2}).
\end{equation*}
This completes the proof of Proposition~\ref{Main-Prop} and hence
of Theorem~\ref{T:Main} of the introduction. The main task is,
thus, to establish the sum formulas (\ref{Sum-1}), (\ref{Sum-2}),
(\ref{Sum-1-P}), and (\ref{Sum-2-P}).

We end this section with the proof of (\ref{rel-qmod}).

\begin{lemma}
\label{L:qmf}\
$4\,t^{2}\,G_{2}^{\prime}(t^{2})\,=\,32\,G_{2}^{2}(t^{2})\,-\,
40\,G_{2}(t)\,G_{2}(t^{2})\,+\,8\,G_{2}^{2}(t)\,-\,t\,G_{2}^{\prime}(t)$.
\end{lemma}

\begin{proof} Let $G_{2}(z)=\sum_{d\geq 0}\sigma(d)\,q^{d}$
and set
\begin{equation*}
 E(z)\,=\,-2\,DG_{2}(2z)\,+\, 32\,G_{2}^{2}(2z)\,-\,
 40\,G_{2}(z)\,G_{2}(2z)\, +\,8\,G_{2}^{2}(z)\,-\,DG_{2}(z)
\end{equation*}
where $z\in {\Bbb C}$ with $\mbox{Im(z)}>0$, $q=e^{2\pi i z}$ and
$D=q\frac{d}{dq}$ is the logarithmic differential operator.
It then suffices to show that $E(z)\equiv 0$.

Since the Eisenstein series $G_{2}(z)$ of weight 2 satisfies
\begin{equation*}\label{G-2rel}
  G_{2}\Big(\frac{az+b}{cz+d}\Big)\,=\,
 (cz+d)^{2}\,G_{2}(z)\,-\, \frac{c(cz+d)}{4\pi i}
 \ \ \ \mbox{for\  \,any\ }\
 \Big(
  \begin{array}{rr}
    a & b \\
    c & d
  \end{array}
 \Big)\,\in\,SL(2,{\Bbb Z})
\end{equation*}
one can show by hand that $E(z)$ is a modular form of
weight $4$ and of level $2$ on the Hecke subgroup $\Gamma_{0}(2)$.
The space of such modular forms is a 2-dimensional vector space
with  generators
\begin{align*}
 G_{4}(z)\,&=\,\frac{1}{24}\,+\,10\sum_{d\geq 1}\,\sigma_{3}(d)\,q^{d}
 \,=\,\frac{1}{24}\,+\,10\,q\,+\,90\,q^{2}\,+\,\cdots
 \\
 G_{2}^{(4)}(z)\,&=\,\big[ G_{2}(z)\,-\,2\,G_{2}(2z)\,\big]^{2}\,=\,
 \frac{1}{24^{2}}\,+\,\frac{1}{12}\,q\,+\,
 \frac{26}{24}\,q^{2}\,+\,\cdots
\end{align*}
where $\sigma_{3}(d)=\sum_{k|d}k^{3}$\ (cf \cite{k}).
Thus, $E(z)$ can be written as
\begin{equation*}
 E(z)\,=\,a\,G_{4}(z)\,+\,b\,G_{2}^{(4)}(z)\,=\,
 \Big(\frac{a}{24}\,+\,\frac{b}{24^{2}}\Big)\,+\,
 \Big(10\,a\,+\,\frac{b}{12}\Big)\,q\,+\,\cdots
\end{equation*}
for some constants $a$ and $b$. On the other hand,
one can also show by hand that the first two terms
in the $q$--expansion of $E(z)$  vanish. Consequently, $a=b=0$ and hence
$E(z)\equiv 0$.
\end{proof}


\setcounter{equation}{0}
\section{Symplectic sum formula}
\label{section3}

Let $E(0)=S^{2}\times T^{2}\to S^{2}$ be a rational elliptic
surface. To save notation, we also denote by $S$ and $F$  the
section class and the fiber class of $E(0)$, respectively.  In
this section, we apply the sum formula \cite {ip3} of the
symplectic sum of $E(2)$ and $E(0)$  to prove the sum formulas
(\ref {Sum-1}), (\ref{Sum-2}), (\ref{Sum-1-P}), and
(\ref{Sum-2-P}).

Recall that $V\cong T^{2}$ is a fixed fiber of $E(2)\to S^{2}$.
For convenience, we use the same notation $V$ for a fixed fiber of
$E(0)\to S^{2}$.  Recall that relative GT (Gromov--Taubes)
invariants of $(E(0),V)$ count $V$--regular maps from possibly
disconnected domain \cite{ip2}. The sum formula of \cite{ip3}
applies to GT invariants to give relations between GT invariants
and relative GT invariants.
The sum formula also applies to GW invariants. In this case, it
 gives relations between GW invariants and (partial)
relative GT invariants --- these invariants are defined to  count
maps each of whose domain component has contact order at least one
with $V$. We denote such (partial) relative invariants of
$(E(0),V)$ for the class $2S+dF$ with the Euler characteristic
$\chi$ and the multiplicity vector $s$ by
\begin{equation}\label{GTI}
 G\Phi^{V}_{2S+dF,\chi,s}\big(C_{\g_{1}\cdots\g_{l}};\kappa;
                              \b_{1},\cdots,\b_{k}\big).
\end{equation}
Here, $s=(s_1,\cdots,s_l)$ equals (1,1) or (2), $\g_{j}\in H_{*}(V;{\Bbb
Z})$, $\kappa\in H_{*}(\CM_{\chi,k+l};{\Bbb Q})$ and $\b_{i}\in
H_{*}(E(0);{\Bbb Z})$; $\CM_{\chi,k+l}$ is the space of all compact
Riemann surface of Euler characteristic $\chi$ with $k+l$ marked
points. We will also use the notation $\g_{j}^{n}$ if the constraint $\g_{j}$
repeats $n$ times and $\mbox{deg}(\g_{j}^{*})$ is even. By
dimension formula of \cite{ip2}, the (partial) GT invariant
(\ref{GTI}) vanishes unless
\begin{equation*}
 \mbox{deg}(\kappa^{*})\,+\,\sum\mbox{deg}(\b_{i}^{*})
 \,+\,\sum\mbox{deg}(\g_{j}^{*})\,=\,
 2\big(4-\tfrac{1}{2}\chi+k+l-2\big).
\end{equation*}

Consider the symplectic sum of $E(2)$ and $E(0)$ along $V$
\begin{equation}  \label{symp-sum}
 E(2)\,=\,E(2)\,\,\#_{V}\,E(0).
\end{equation}
The Gluing Theorem (Theorem 10.1 of \cite{ip3}) applies for this
sum to  give relations between family GW invariants of $E(2)$ for
the classes $2S+dF$ and (partial) relative family GT invariants of
$(E(2),V)$ for the classes $2S+dF$. The latter count $V$--regular
maps $f$ with possibly disconnected domains such that if $f$ has a
disconnected domain of Euler characteristic $\chi$ with $k$ marked
points $f$ is a pair of $V$--regular maps $(f_{1},f_{2})$
satisfying $[f_{i}]=S+d_{i}F$ with $d_{1}+d_{2}=d$ and the domain
of $f_{i}$ lies in $\CM_{g_{i},k_{i}}$ with
$\chi=4-2(g_{1}+g_{2})$ and $k=k_{1}+k_{2}$. Denote the moduli
space of all such pairs $(f_{1},f_{2})$ by
\begin{equation}\label{cor-moduli}
 \CM_{(g_{1},g_{2}),(k_{1},k_{2}),(1,1)}^{\H,V}\big(d_{1},d_{2}\big).
\end{equation}
The standard corbodism argument (cf proof of Proposition 3.7 of
\cite{l1}) then shows that the moduli space (\ref{cor-moduli}) is
corbodant to the product moduli space
$$ \CM_{g_{1},k_{1},(1)}^{\H,V}(S+d_{1}F)\,\times \,
   \CM_{g_{2},k_{2},(1)}^{V}(S+d_{2}F)  $$
where the first factor is a relative family GW moduli space of
$(E(2),V)$ and the second is a relative ordinary GW moduli space
of $(E(2),V)$. Since $E(2)$ is a K3 surface, the contribution of
the second factor to relative ordinary GW invariants of $(E(2),V)$
vanishes. Consequently,  the contribution of the moduli space
(\ref{cor-moduli}) to (partial) relative family GT invariants also vanishes.
This implies  that for the classes $2S+dF$
the (partial) relative family GT invariants of $(E(2),V)$
and the relative family GW invariants $GW^V$ are the same.
Therefore, a family version of the sum formula of the sum (\ref{symp-sum})
for the classes $2S+dF$ relates
$GW^{\H}$ invariants of $E(2)$ and $GW^{V}$ invariants  of $(E(2),V)$.

On the other hand, using Lemma 14.5 of \cite{ip3} and routine
dimension count one can show that there is no `contribution from
the neck' (cf Section 12 of \cite{ip3}). Moreover, `rim tori'
\cite{ip2} of $E(2)$--side disappear under the symplectic sum (\ref
{symp-sum}) --- that enables us to work with summed relative
invariants. Combined with these observations, the Gluing Theorem
then yields a considerably simple sum formulas (\ref{fv-sum-for})
below: Let $\{\g_{i}\}$ be a basis of $H_{*}(V;{\Bbb Z})$ and
$\{\g^{i}\}$ be its dual basis with respect to the intersection
form of $V$. For a vector of nonnegative integers
$m=(m_{1},\cdots,m_{4})$   with $\sum\,m_{i}$ is either 2 or 1, we
set
\begin{equation*}
 C_{\g_{m}}\,=\,C_{\g_{1}^{m_{1}}\cdots\g_{4}^{m_{4}}}\,,\ \
 C_{\g_{m^{*}}}=\,C_{(\g^{4})^{m_{4}}\cdots(\g^{1})^{m_{1}}}\,,\ \
 \mbox{and}\ \
 m!\,=\,\prod\,m_{i}!\,.
\end{equation*}
For a multiplicity vector $s=(s_{1},\cdots,s_{l})$, either (1,1)
or (2), let $|s|=\prod s_{i}$. We are now ready to write a sum
formula of the symplectic sum (\ref {symp-sum})
\begin{align}  \label{fv-sum-for}
 &GW_{2S+dF,g}^{{\cal H}}
 \big(\tau(F)^k,pt^{g-k}\big)  \nonumber \\
 &=\,\sum\,\frac{|s|}{m!}\,GW^{V}_{2S+d_{1}F,g_{1},s}
 \big( C_{\gamma_{m}}\big)\, G\Phi^{V}_{2S+d_{2}F,\chi_{2},s}
 \big( C_{\gamma_{m^{*}}};\tau(F)^k,pt^{g-k}\big)
\end{align}
where  the sum is over all $s=(s_{1},\cdots,s_{l})$ which is
either (1,1) or (2), vectors $m=(m_{i})$ as above with $\sum
m_{i}=l(s)$, $d=d_{1}+d_{2}$ and
$g=g_{1}-\frac{1}{2}\chi_{2}+l(s)$.

Similarly, one can also derive a sum formula for the case of
family invariants of $E(2)$ for the classes $S+(2d-3)F$, $d\in
{\Bbb Z}$. Recall that $U\cong T^{2}$ is a fixed symplectic
submanifold of $E(2)$ that represents the class $2F$. We also
denote by $U$ a fixed fiber of $E(0)\to S^{2}$ and
 consider the symplectic sum of $E(2)$ with $E(0)$ along $U$
\begin{equation}\label{SS-U}
 E(2)\,=\,E(2)\,\,\#_{U}\,E(0).
\end{equation}
Since both $U$ and $V$ are fibers of $E(0)$, relative GW invariants
of $(E(0),V)$ and $(E(0),U)$ are in fact the same. Thus, the
(partial) relative  GT invariants of $(E(0),V)$ and $(E(0),U)$ are
also the same invariants, ie $G\Phi^{V}=G\Phi^{U}$. On the other
hand, $U\subset E(2)$ represents $2F$ on $E(2)$, while $U\subset
E(0)$ represents $F$ on $E(0)$. With these observations, repeating
the same arguments as above for the symplectic sum (\ref{SS-U})
gives a sum formula like (\ref{fv-sum-for})
\begin{align}  \label{fv-sum-for-U}
  & GW_{S+(2d-3)F,g}^{{\cal H}}
  \big(\tau(2F)^k,pt^{g-k}\big)  \nonumber \\
  &=\,
 \sum\,\frac{|s|}{m!}\,GW^{U}_{S+(2d_{1}-3)F,g_{1},s}
 \big( C_{\gamma_{m}}\big)\,G\Phi^{V}_{2S+d_{2}F,\chi_{2},s}
 \big( C_{\gamma_{m^{*}}};\tau(F)^k,pt^{g-k}\big).
\end{align}

\begin{rem}\label{sd}
Once and for all, we fix  an (ordered) basis $\{\,
pt,\g_{1},\g_{2},F\, \}$ of $H_{*}(V;{\Bbb Z})\cong H_{*}(U;{\Bbb
Z})$ and its (ordered) dual basis $\{\, F,\g_{2},-\g_{1},pt\, \}$
with respect to the intersection form of $V\cong U$ where
$\{ \g_{1}, \g_{2} \}$ is a basis of
 $H_{1}\big(V;{\Bbb Z} \big)\cong H_{1}\big(U;{\Bbb Z} \big)\cong
 H_{1}\big(E(0);{\Bbb Z}\big)$ with $\g_{1}\cdot \g_{2}=1$.
Then,
in the sum formulas (\ref{fv-sum-for}) and (\ref{fv-sum-for-U})
the splitting of diagonal for contact constraints $C_{\g_m}$
is given as follows:
\begin{enumerate}
\item[$\bullet$] if $m=(2,0,0,0)$ then $\g_m=pt^2$  and $\g_{m^*}=F^2$,
\item[$\bullet$] if $m=(1,0,0,1)$ then $\g_m=pt\cdot F$ and $\g_{m^*}=pt\cdot F$,
\item[$\bullet$] if $m=(0,1,1,0)$ then $\g_{m}=\g_{1}\cdot\g_{2}$ and
                 $\g_{m^*}=(-\g_1)\cdot \g_2$.
\end{enumerate}
\end{rem}

Using the sum formula (\ref{fv-sum-for}), one can derive relations
between invariants $GW^{\H}$ and $GW^{V}$.

\begin{lemma}
\label{gw-rel(gw)} Let $GW^{V}$ be the  relative family invariants
of $(E(2),V)$. Then,

\begin{enumerate}
\item[\rm(a)]  $\frac{1}{2}\,GW^{V}_{2S+dF,0,(1,1)}\big(C_{F^{2}}\big)\,=\,
             GW^{\H}_{2S+dF,0}$\,,
\item[\rm(b)]  $GW^{V}_{2S+dF,1,(1,1)}\big(C_{\g_{1}\cdot \g_{2}}\big)\,=\,
            GW^{V}_{2S+dF,1,(1,1)}\big(C_{pt\cdot F}\big)$\,,
\item[\rm(c)]  $GW^{V}_{2S+dF,1,(1,1)}\big(C_{pt\cdot F}\big)\,=\,
            GW^{\H}_{2S+dF,1}\big(pt\big)\,-\,
            2\!\!\underset{d=d_1 +d_2}{\sum}\!
            GW^{\H}_{2S+d_{1}F,0}\,d_{2}\,\sigma(d_{2})$.
\end{enumerate}
\end{lemma}

\begin{proof}
(a)\qua By the sum formula (\ref{fv-sum-for}), Remark~\ref{sd} and
Lemma~\ref{RI-L3.1}\,a, we have
\begin{align*}
 GW^{{\cal H}}_{2S+dF,0}\,&=\,
 \sum_{d=d_{1}+d_{2}}
 \tfrac{1}{2}\, GW^{V}_{2S+d_{1}F,0,(1,1)}\big(C_{F^{2}}
 \big)\,G\Phi^{V}_{2S+d_{2}F,4,(1,1)}\big(C_{pt^{2}}\big) \\
 &=\,
 \tfrac{1}{2}\, GW^{V}_{2S+dF,0,(1,1)}\big(C_{F^{2}}\big).
\end{align*}

\medskip
(b)\qua Note that the sum formula (\ref{fv-sum-for}) also holds for
one dimensional constraints. The sum formula (\ref{fv-sum-for}),
Remark~\ref{sd}
and Lemma~\ref{RI-L3.1}\,b,c,d,e thus give
\begin{align*}
 &GW^{\H}_{2S+dF,1}\big(\iota_*(\g_{1}),\iota_*(\g_{2})\big)
 \nonumber \\
  &= \sum_{d=d_{1}+d_{2}}
    GW^{V}_{2S+d_{1}F,1,(1,1)}\big(C_{pt\cdot F}\big)
    \, G\Phi^{V}_{2S+d_{2}F,4,(1,1)}\big(C_{pt\cdot F};\g_{1},\g_{2}\big)
 \nonumber \\
 &+ \sum_{d=d_{1}+d_{2}}
   GW^{V}_{2S+d_{1}F,1,(1,1)}\big(C_{\g_{1}\cdot \g_{2}} \big)\,
   G\Phi^{V}_{2S+d_{2}F,4,(1,1)}
       \big(C_{(-\g_{1})\cdot \g_{2}};\g_{1},\g_{2}\big)
 \nonumber \\
 &+ \sum_{d=d_{1}+d_{2}}2\,
   GW^{V}_{2S+d_{1}F,1,(2)}\big(C_{F}\big)\,
   G\Phi^{V}_{2S+d_{2}F,2,(2)}\big(C_{pt};\g_{1},\g_{2}\big)
 \nonumber \\
 &+\,\sum_{d=d_{1}+d_{2}}\tfrac{1}{2}\,
   GW^{V}_{2S+d_{1}F,0,(1,1)}\big( C_{F^{2}}\big)\,
   G\Phi^{V}_{2S+d_{2}F,2,(1,1)}\big(C_{pt^{2}};\g_{1},\g_{2} \big)
 \nonumber \\
  &=\ \ GW^{V}_{2S+dF,1,(1,1)}\big(C_{pt\cdot F}\big)\,-\,
        GW^{V}_{2S+dF,1,(1,1)}\big(C_{\g_{1}\cdot \g_{2}} \big)
\end{align*}
where $\iota\co V\hookrightarrow E(2)$ is the inclusion map. Since
$E(2)$ is simply connected, the left hand side of the first
equality vanishes and hence this shows (b).

\medskip

(c)\qua We have
\begin{align*}
 \hspace{-0.5cm}
 &GW^{\H}_{2S+dF,1}\big(pt\big)\\
 &=\, \sum_{d=d_{1}+d_{2}}
   GW^{V}_{2S+d_{1}F,1,(1,1)}\big(C_{pt\cdot F}\big)\,
   G\Phi^{V}_{2S+d_{2}F,4,(1,1)}\big(C_{pt\cdot F};pt\big)
 \nonumber \\
 &+\, \sum_{d=d_{1}+d_{2}}2\,
   GW^{V}_{2S+d_{1}F,1,(2)}\big(C_{F}\big)\,
   G\Phi^{V}_{2S+d_{2}F,2,(2)}\big(C_{pt};pt\big)
 \nonumber \\
 &+\, \sum_{d=d_{1}+d_{2}}\tfrac{1}{2}\,
   GW^{V}_{2S+d_{1}F,0,(1,1)}\big( C_{F^{2}}\big)\,
   G\Phi^{V}_{2S+d_{2}F,2,(1,1)}\big(C_{pt^{2}};pt\big)
 \notag \\
  &=\,
   GW^{V}_{2S+dF,1,(1,1)}\big(C_{pt\cdot F}\big)\,+\,
   2\underset{d=d_{1}+d_{2}}{\sum}
   GW^{\H}_{2S+d_{1}F,0}\, d_{2}\,\sigma(d_{2})
\end{align*}
where the first equality follows from the sum formula (\ref{fv-sum-for})
and Remark~\ref{sd},
and the second follows from Lemma~\ref{RI-L3.1}\,f,g,h.
\end{proof}

The sum formula (\ref{fv-sum-for-U}) also gives relations
between invariants $GW^{\H}$ and $GW^{U}$.

\begin{lemma}
\label{gw-rel(gw)-U} Let $GW^{U}$ be the  relative family
invariants of $(E(2),U)$. Then,
\begin{enumerate}
\item[\rm(a)]  $\frac{1}{2}\,GW^{U}_{S+(2d-3)F,0,(1,1)}\big(C_{F^{2}}\big)\,=\,
             GW^{\H}_{S+(2d-3)F,0}$\,,
\item[\rm(b)]  $GW^{U}_{S+(2d-3)F,1,(1,1)}\big(C_{\g_{1}\cdot \g_{2}}\big)\,=\,
            GW^{U}_{S+(2d-3)F,1,(1,1)}\big(C_{pt\cdot F}\big)$\,,
\item[\rm(c)]  $ GW^{U}_{S+(2d-3)F,1,(1,1)}\big(C_{pt\cdot F}\big)$
\item[]     $=\, GW^{\H}_{S+(2d-3)F,1}\big(pt\big)\,-\,2
             \underset{d=d_{1}+d_2}{\sum}
            GW^{\H}_{S+(2d_1 -3) F,0}\,d_{2}\,\sigma(d_{2})$.
\end{enumerate}
\end{lemma}

The proof of this lemma is identical to the proof of
Lemma~\ref{gw-rel(gw)}.

Now, we are ready to show the sum formulas (\ref{Sum-1})
and (\ref{Sum-2}).

\begin{prop}
\label{P:sum-formula} Let $M$ and $M^V$ be the generating
functions defined in (\ref{def-M}) and (\ref{def-M^V}),
respectively. Then,
\begin{enumerate}
\item[\rm(a)]  $M_{1}\big(\tau(F)\big)\,=\,
             M_{1,(2)}^{V}\, +\, 4\,M_{0}\,G_{2}$,
\item[\rm(b)]  $M_{2}\big(\tau(F)^2\big)\,=
             \,2\,M_{1}(pt)\,+\,20\,M^{V}_{1,(2)}\,G_{2}\,+\,
             M_{0}\left(16\,G_{2}^{2}\,+\,8\,t\,G_{2}^{\prime}\right).$
\end{enumerate}
\end{prop}
\begin{proof} (a)\qua We have
\begin{align}\label{P3.3a}
 &GW^{\H}_{2S+dF,1}\big(\tau(F)\big)
 \notag \\
 &=\  \sum_{d=d_{1}+d_{2}}
 GW^{V}_{2S+d_{1}F,1,(1,1)}\big(C_{pt\cdot F}\big)\,
 G\Phi^{V}_{2S+d_{2}F,4,(1,1)} \big(C_{pt\cdot F};\tau(F)\big)
 \nonumber \\
 &+\ \sum_{d=d_{1}+d_{2}} 2\,
 GW^{V}_{2S+d_{1}F,1,(2)}\big(C_{F}\big)\,
 G\Phi^{V}_{2S+d_{2}F,0,(2)} \big(C_{pt};\tau(F)\big)
 \nonumber \\
 &+\ \sum_{d=d_{1}+d_{2}}
 \tfrac{1}{2}\, GW^{V}_{2S+d_{1}F,0,(1,1)}\big(C_{F^{2}} \big)\,
 G\Phi^{V}_{2S+d_{2}F,2,(1,1)} \big(C_{pt^{2}};\tau(F)\big).
 \notag \\
 &=\ GW^{V}_{2S+dF,1,(2)}\big(C_{F}\big)\,+\,
 \sum_{d=d_{1}+d_{2}}2\,
 GW^{V}_{2S+d_{1}F,0,(1,1)}\big(C_{F^{2}} \big)\,\sigma(d_{2})
 \notag \\
 &=\ GW^{V}_{2S+dF,1,(2)}\big(C_{F}\big)\,+\,
 \sum_{d=d_{1}+d_{2}}4\, GW^{\H}_{2S+d_{1}F,0}\,\sigma(d_{2})
 \end{align}
where the first equality follows from the sum formula
(\ref{fv-sum-for}) and Remark~\ref{sd}, the second equality follows from
Lemma~\ref{RI-P3.3a} and the third equality follows from
Lemma~\ref{gw-rel(gw)}\,a. Then, (a) follows from (\ref{P3.3a})
and definition of generating functions.

\medskip 
(b)\qua The sum formula (\ref{fv-sum-for}) and Remark~\ref{sd} give
\begin{align}  \label{pf-main-a3}
 \hspace{-0.2cm}
 &GW^{\H}_{2S+dF,2}\big(\tau(F)^2\big)
 \notag \\
 &=\sum_{d=d_{1}+d_{2}}
 \tfrac{1}{2}\,GW^{V}_{2S+d_{1}F,2,(1,1)}\big(C_{pt^{2}} \big)\,
 G\Phi^{V}_{2S+d_{2}F,4,(1,1)} \big(C_{F^{2}};\tau(F)^2 \big)
  \notag \\
 &+\sum_{d=d_{1}+d_{2}}
 2\,GW^{V}_{2S+d_{1}F,2,(2)}\big(C_{pt}\big)\,
 G\Phi^{V}_{2S+d_{2}F,0,(2)} \big(C_{F};\tau(F)^2\big)
 \notag\\
 &+\sum_{d=d_{1}+d_{2}}
 GW^{V}_{2S+d_{1}F,1,(1,1)}\big(C_{pt\cdot F}\big)\,
 G\Phi^{V}_{2S+d_{2}F,2,(1,1)} \big(C_{pt\cdot F};\tau(F)^2\big)
 \notag \\
  &+\sum_{d=d_{1}+d_{2}}
 GW^{V}_{2S+d_{1}F,1,(1,1)}\big(C_{\g_{1}\cdot \g_{2}} \big)\,
 G\Phi^{V}_{2S+d_{2}F,2,(1,1)}
    \big(C_{(-\g_{1})\cdot  \g_{2}};\tau(F)^2\big)
  \notag \\
 &+\sum_{d=d_{1}+d_{2}}
 2\,GW^{V}_{2S+d_{1}F,1,(2)}\big(C_{F}\big)\,
 G\Phi^{V}_{2S+d_{2}F,1,(2)} \big(C_{pt};\tau(F)^2\big)
 \notag \\
 &+\sum_{d=d_{1}+d_{2}}
 \tfrac{1}{2}\, GW^{V}_{2S+d_{1}F,0,(1,1)}\big(C_{F^{2}} \big)\,
 G\Phi^{V}_{2S+d_{2}F,0,(1,1)}
   \big(C_{pt^{2}};\tau(F)^2 \big)
 \end{align}
By Lemma~\ref{RI-P3.3b}, the right hand side of (\ref{pf-main-a3})
becomes
\begin{align*}
 &GW^{V}_{2S+dF,1,(1,1)}\big(C_{pt\cdot F}\big)\,+\,
 GW^{V}_{2S+dF,1,(1,1)}\big(C_{\g_{1}\cdot \g_{2}} \big)\\
 +&
 \sum_{d=d_{1}+d_{2}}
 20\,GW^{V}_{2S+d_{1}F,1,(2)}\big(C_{F}\big)\,\sigma(d_{2})
 \notag \\
 +&\sum_{d=d_{1}+d_{2}}
 \tfrac{1}{2}\, GW^{V}_{2S+d_{1}F,0,(1,1)}\big(C_{F^{2}} \big)\,
 \Big(  \sum_{k_{1}+k_{2}=d_{2}}16\,\sigma(k_{1})\,\sigma(k_{2})\,
  +\,12\,d_{2}\,\sigma(d_{2})\, \Big).
\end{align*}
This can be further simplified by using Lemma~\ref{gw-rel(gw)} to
give
\begin{align}\label{P3.3-bl}
 &2\,GW^{\H}_{2S+dF,1}\big(pt\big)\,+
 \sum_{d=d_{1}+d_{2}}
 20\,GW^{V}_{2S+d_{1}F,1,(2)}\big(C_{F}\big)\cdot
 \sigma(d_{2})  \notag \\
 &+\,\sum_{d=d_{1}+d_{2}}
 GW^{\H}_{2S+d_{1}F,0}
 \Big(  \sum_{k_{1}+k_{2}=d_{2}}16\,\sigma(k_{1})\,\sigma(k_{2})\,
  +\,8\,d_{2}\,\sigma(d_{2})\, \Big).
\end{align}
Thus, (b) follows from  (\ref{pf-main-a3}), (\ref{P3.3-bl}) and
definition of generating functions.
\end{proof}

The same computation shows (\ref{Sum-1-P}) and (\ref{Sum-2-P}):

\begin{prop}
\label{P:sum-formula-U} Let $P$ and $P^U$ be the generating
functions defined in (\ref{def-P}) and (\ref{def-P^U}),
respectively. Then,

\begin{enumerate}
\item[\rm(a)]  $P_{1}\big(\tau(2F)\big)\,=\,
             P_{1,(2)}^{U}\, +\, 4\,P_{0}\,G_{2}$,
\item[\rm(b)]  $P_{2}\big(\tau(2F)^2\big)\,=
             \,2\,P_{1}(pt)\,+\,20\,P^{U}_{1,(2)}\,G_{2}\,+\,
             P_{0}\left(16\,G_{2}^{2}\,+\,8\,t\,G_{2}^{\prime}\right).$
\end{enumerate}
\end{prop}


\setcounter{equation}{0}
\section{GW invariants of $E(0)$}
\label{section4}

In order to complete the proof of the sum formulas (\ref{Sum-1}),
(\ref{Sum-2}), (\ref{Sum-1-P}), and (\ref{Sum-2-P}), we need to
compute the (partial) Gromov--Taubes invariants of $E(0)$ that appeared in
Section \ref{section3}. Those invariants are expressed in terms of the relative
invariants of $E(0)$. The aim of this section is to compute
various GW invariants of $E(0)$ which we use in later sections to
compute the required relative invariants.

Recall that $S$ and $F$ denote the section class and the fiber
class of $E(0)$, respectively. We always denote  the genus $g$ GW
invariants of $E(0)$ for the class $A$ by
\begin{equation*}
 \Phi_{A,g}\big(\kappa;\b_{1},\cdots,\b_{k}\big)
\end{equation*}
where $\kappa\in H_{*}(\CM_{g,k};{\Bbb Q})$ and $\b_{i}\in
H_{*}(E(0);{\Bbb Z})$. Note that by dimension count this invariant
vanishes unless
\begin{equation*}
 \mbox{deg}(\kappa^{*})\,+\,\sum\mbox{deg}(\b_{i}^{*})\,=\,
 2\,(\,A\cdot (2F)+g-1+k\,)
\end{equation*}
where $\kappa^{*}$ and $\b_{i}^{*}$ are Poincar\'{e} dual of
$\kappa$ and $\b_{i}$, respectively.

We start with the genus 0 GW invariants for the trivial homology
class. The lemma below directly follows from Proposition 1.2 of
\cite{km1}.

\begin{lemma}\label{GW=0}
Let $\Phi$ denote the GW invariants of $E(0)$. Then,
\begin{equation*}
 \Phi_{0,0}\big(\kappa;\b_{1},\cdots,\b_{n}\big) = 0
 \   \mbox{unless}\   n=3\ \ \
\mbox{and}\ \ \
 \Phi_{0,0}\big(\b_{1},\b_{2},\b_{3}\big) =
 \int_{E(0)}\b_{1}^{*}\cup\b_{2}^{*}\cup\b_{3}^{*}
\end{equation*}
where $\b_{i}^{*}$ denote the Poincar\'{e} dual of $\b_{i}\in
H_{*}(E(0);{\Bbb Z})$.
\end{lemma}

Recall that $V\cong T^{2}$ is a fixed fiber of $E(0)\to S^{2}$. We
always denote the genus $g$ relative GW invariants of $(E(0),V)$
for the class $A$ with the multiplicity vector $s$ by
\begin{equation*}
 \Phi^{V}_{A,g,s}\big(\kappa;\b_{1},\cdots,\b_{k};C_{\g_{1}\cdots\g_{l}}\big)
\end{equation*}
where $s=(s_{1},\cdots,s_{l})$ with $\sum s_{j}=A\cdot [V]=A\cdot
F$ and $\g_{j}\in H_{*}(V;{\Bbb Z})$. By dimension formula of
\cite{ip2}, this relative invariant vanishes unless
\begin{equation*}
  \mbox{deg}(\kappa^{*})\,+\,\sum\mbox{deg}(\b_{i}^{*})\,+\,
  \sum\mbox{deg}(\g_{j}^{*})\,=\,
 2\,\big(\,A\cdot (2F)+g-1+k+l-A\cdot F\,\big)
\end{equation*}
where $\g_{i}^{*}$ is the Poincar\'{e} dual of $\g_{i}$.

Recall that  $\{\g_{1},\g_{2}\}$ is a basis of
 $H_{1}(V;{\Bbb Z})\cong H_{1}(E(0);{\Bbb Z})$ with $\g_1\cdot\g_2=1$
and $F$ also denotes the fundamental class  of $V$.

\begin{lemma}{\rm\cite{ip3,ll,l2}}\qua
Let $\Phi$ and $\Phi^V$ denote the GW invariants of $E(0)$ and the
relative GW invariants of $(E(0),V)$, respectively. Then,
\label{well-known}
\begin{enumerate}
\item[\rm(a)]  $\Phi_{S,0}\big(pt\big)\,=\,
                   \Phi_{S,0,(1)}^{V}\big(C_{pt}\big)\,=\,
                   \Phi_{S,0,(1)}^{V}\big(pt;C_{F}\big)\,=\,1$,

\item[\rm(b)]  $\Phi_{S,0}\big(\g_{1},\g_{2}\big)\,=\,
                   \Phi_{S,0,(1)}^{V}\big(\g_{1};C_{\g_{2}}\big)\,=\,
                   \Phi_{S,0,(1)}^{V}\big(\g_{1},\g_{2};C_{F}\big)\,=\,1$,
\item[\rm(c)]  $\Phi_{S+dF,1}\big(\tau(F),pt\big)\,=\,
                   \Phi^{V}_{S+dF,1,(1)}\big(\tau(F);C_{pt}\big)$
\item[]             $=\,
                   \Phi^{V}_{S+dF,1,(1)}\big(\tau(F),pt;C_{F}\big)\,=\,
                   2\,\sigma(d)$,
\item[\rm(d)]  $\Phi_{S+dF,1}\big(pt^{2}\big)=\,2\,d\,\sigma(d)$
                  \,{\rm ;}\,\
                  $\Phi^{V}_{S+dF,1}\big(pt;C_{pt}\big)=\,
                  d\,\sigma(d)$,

\item[\rm(e)]  $\Phi_{S+dF,1}\big(pt,\g_{1},\g_{2}\big)\,=\,
                   \Phi^{V}_{S+dF,1,(1)}\big(pt,\g_{1};C_{\g_{2}}\big)
                   \,=\,d\,\sigma(d)\, {\rm ;}$
\item[]     $\Phi^{V}_{S+dF,1,(1)}\big(\g_{1},\g_{2};C_{pt}\big)=\, 0$,

\item[\rm(f)]  $\Phi_{dF,1}\big(S\big)\,=\,2\,\sigma(d)$.
\end{enumerate}
\end{lemma}

We will often use the following simple observations:

\begin{rem}
\label{rat-cur-E(0)} Consider $E(0)=S^{2}\times T^{2}$ with a
product complex structure. Since there is no nontrivial
holomorphic map from $S^{2}$ to $T^{2}$, any nontrivial
holomorphic map $f\co S^{2}\to E(0)$ should represent a class $aS$,
$a\geq 1$, with its image a section. Thus, the image of such maps
can't pass through generic two points, a generic geometric
representative of $\g_{1}$ or $\g_{2}$ and a generic point, or a
generic section and a generic point. Combined with the Gromov
Convergence Theorem \cite{pw,p,is},\ this observation gives
vanishing results of certain genus 0 invariants. For example,
$\Phi_{aS+dF,0}\big(\,\cdot\,\big)\,
 =\, \delta_{d0}\,\Phi_{aS,0}\big(\,\cdot\,\big)$ and
\begin{equation*}
 \Phi_{aS,0}\big(pt,\g_{1},\,\cdot\,\,\big)\,=\,
 \Phi_{aS,0}\big(S,pt,\,\cdot\,\,\big)\,=\,
 \Phi_{aS,0}\big(S,S,\,\cdot\,\,\big)\,=\,0\,.
\end{equation*}
\end{rem}

Fix a product complex structure on $E(0)$ and let
 $f\co S^2\to E(0)$ be a  holomorphic map  representing a class $aS$, $a\geq 1$.
Then $f$ is a branched covering of some section $S_0$ of $E(0)$
and since the normal
bundle of $S_0$ is trivial $H^1(f^*TE(0))=H^1(\co\oplus f^*TS_0)=0$.
This shows that the linearization
$L_f$ of the holomorphic map equation at $f$ has a trivial
cokernal. With this observation, we will use the product complex
structure of $E(0)$ for computation of both absolute and relative
GW invariants for the following cases:

\begin{lemma}
\label{lemma5-1} Let $\Phi$ and $\Phi^V$ denote the GW invariants
of $E(0)$ and the relative GW invariants of $(E(0),V)$,
respectively. Then, we have

\begin{enumerate}
\item[\rm(a)]  $2\,\Phi^{V}_{2S,0,(2)}\big(\tau(F);C_{pt}\big)\,=\,
                  \Phi^{V}_{2S,0,(2)}\big(pt,\tau(F);C_{F}\big)\,=\,1$,

\item[\rm(b)]  $\Phi^{V}_{S,0,(1)}\big(\tau(F);C_{F}\big)\,=\,
                   \Phi_{S,0}\big(\tau(F))\,=\,
                  \Phi^{V}_{2S,0,(2)}
                    \big(\tau(F)^2;C_{F}\big)\,=\,0$.
\end{enumerate}
\end{lemma}

\begin{proof} (a)\qua Fix a product complex structure on
$E(0)$ and let $V_{1},V_{2},V$ be distinct fibers and $p$ be a
point in $V$. Denote by
\begin{equation}  \label{cut-down1}
 \M^{V}_{0,2,(2)}\big(E(0),2S\big)(V_1,V_2,p)\ \subset \
 \CM_{0,3}\big(E(0),2S\big)
\end{equation}
the cut-down moduli space that consists of all maps
$(f,C;x_{1},x_{2},x_{3})$
satisfying (i) the contact order of $f$ with $V$ at $x_{3}$ is 2
and (ii) $ f(x_{i})\in V_{i}$, $i=1,2$, and $f(x_{3})=p$. This
space is smooth of expected complex dimension 1 with no boundary
stratum; by stability each map in the space (\ref{cut-down1}) has
a smooth domain. Moreover, each map $f$ in the space
(\ref{cut-down1}) has no non-trivial automorphism and the
linearization $L_f$ has a trivial cokernal. Therefore, the
invariant $\Phi^{V}_{2S,0,(2)} \big(\tau(F),F;C_{pt}\big)$ is
equal to the (homology) Euler class of the relative cotangent bundle $ {\cal
L}_{1}$ over the  cut-down moduli space (\ref{cut-down1}) whose
fiber over $(f,C;x_{1},x_{2},x_{3})$ is $T^*_{x_1}C$.

Note that each map $f$ in (\ref{cut-down1}) is a 2-fold branched
covering of the fixed section $S_{p}\cong {\Bbb P}^{1}$ containing
the point $p$. Let $p_{1},p_{2}$, and $p_{3}$ be distinct points
of ${\Bbb P}^{1}$ and $P=\{p_{3}\}$. Then the space
(\ref{cut-down1}) can be identified with
\begin{equation}  \label{cut-down2}
 \M^{P}_{0,2,(2)}\big({\Bbb P}^{1},2\big)(p_1,p_2)\
 \subset\ \CM_{0,3}\big({\Bbb P}^1,2\big)
\end{equation}
the space of degree 2 stable maps
 $f\co \big({\Bbb P}^{1},x_{1},x_{2},x_{2}\big) \,\to\,
 \big({\Bbb P}^{1},P\big)$
satisfying (i) the contact order of $f$ with $p_{3}$ at $x_{3}$ is
2, and (ii) $f(x_{i})=p_{i}$ for $i=1,2$. Under this
identification, the relative cotangent bundle ${\cal L}_1$ over
the space (\ref{cut-down1}) becomes the relative cotangent bundle,
still denoted by ${\cal L}_1$, over the space (\ref{cut-down2});
this bundle $ {\cal L}_{1}$ has a fiber $T^{*}_{x_{1}}{\Bbb
P}^{1}$ at $f$.

For each map $f$ in the space (\ref{cut-down2}), choose local
holomorphic coordinates $z$ centered at $x_{1}$ and $w$ centered
at $p_{1}$. Then, there is a local expansion $f(z)\,=\,\sum_{k\geq
1}a_{k}\,z^{k}$. The leading coefficient $a_{1}$ is the 1-jet of
$f$ at $x_{1}$ modulo higher order terms. Thus, we have a global
section
\begin{equation}  \label{g-section}
 a_{1}\,\in\, {\cal L}_{1}
\end{equation}
over the space (\ref{cut-down2}).  The zero set of this section
consists of degree two branched coverings
 $ \big({\Bbb P}^{1},x_{1},x_{2},x_{3}\big)\to
 \big({\Bbb P}^{1},p_{1},p_{2},p_{3}\big)$
with the ramification indexes (2,1,2) at marked points
$(x_{1},x_{2},x_{3})$. Since there is only one such map, the (homology) Euler
class of the bundle ${\cal L}_1$ over the space (\ref{cut-down2})
is one. Consequently, we have
\begin{equation*}
 2\,\Phi^{V}_{2S,0,(2)}\big(\tau(F);C_{pt}\big)\,=\,
 \Phi^{V}_{2S,0,(2)}\big(\tau(F),F;C_{pt}\big)\,=\,1.
\end{equation*}

By the same arguments as above, the invariant
 $\Phi^{V}_{2S,0,(2)}\big(pt,\tau(F);C_{F}\big)$
is the number of degree two branched coverings
 $\big({\Bbb P}^{1},x_{1},x_{2},x_{3}\big)\to
 \big({\Bbb P}^{1},p_{1},p_{2},p_{3}\big)$
with the ramification indexes (1,2,2). Since the number of such
maps is 1, we have
\begin{equation*}
\Phi^{V}_{2S,0,(2)}\big(pt,\tau(F);C_{F}\big)\,=\,1.
\end{equation*}

\medskip
(b)\qua Similarly, as above, one can show that the
invariant $\Phi^{V}_{S,0,(1)}\big(\tau(F),F;C_{F}\big)$ is the
Euler class of the relative cotangent bundle
\begin{equation*}
 {\cal L}_{1}\,\to \,
 \left(\, \M^{P}_{0,2,(1)}\big({\Bbb P}^{1},1\big)(p_1,p_2)\,\right)
 \,\times V
\end{equation*}
with a section defined similarly as in (\ref{g-section}). The zero
set of this section is empty since there is no degree 1 map
   ${\Bbb P}^{1}\to {\Bbb P}^{1}$ with ramification indexes $(2,1,1)$. Thus, we have
\begin{equation*}
 \Phi^{V}_{S,0,(1)}\big(\tau(F);C_{F}\big)\,=\,
 \Phi^{V}_{S,0,(1)}\big(\tau(F),F;C_{F}\big)\,=\,0.
\end{equation*}
Repeating the same argument, we also have
\begin{equation*}
 \Phi_{S,0}\big(\tau(F)\big)\,=\,
 \Phi_{S,0}\big(\tau(F),F^2\big)\,=\,0.
\end{equation*}

Similarly, the invariant
 $\Phi^{V}_{2S,0,(2)}\big(\tau(F)^2;C_{F}\big)$
 is the Euler class of the bundle
\begin{equation*}
 {\cal L}_{1}\,\oplus {\cal L}_{2}\,\to \,\left(\, \M^{P}_{0,2,(2)}
 \big( {\Bbb P}^{1},2\big)(p_1,p_2)\,\right)\,\times V
\end{equation*}
with a section defined similarly as in (\ref{g-section}). The zero
set of this section is empty since there is no degree 2 map
 ${\Bbb P}^{1}\to {\Bbb P}^{1}$ with ramification indexes $(2,2,2)$.
 Thus, the invariant
 $\Phi^{V}_{2S,0,(2)}\big(\tau(F)^2;C_{F}\big)$ is trivial.
\end{proof}


\setcounter{equation}{0}
\section{GW invariants of $E(0)$ with $\tau(F)$ constraints}
\label{section5}

The aim of this section is to prove:

\begin{lemma}\label{L5.1}
Let $\Phi$ denote the GW invariants of $E(0)$. Then,
\begin{enumerate}
\item[\rm(a)] $\Phi_{2S,0}\big(\tau(F)^2,pt\big)\,=\,1$,
\item[\rm(b)] $\Phi_{2S,0}\big(\tau(F)^2,\g_{1},\g_{2}\big)\,=\,2$,
\item[\rm(c)] $\Phi_{2S+dF,1}\big(\tau(F)^3,pt\big) \,=\,
            24\,\sigma(d)$,
\item[\rm(d)] $\Phi_{2S+dF,1}\big(\tau(F)^2,pt^{2}\big)\,=\,
            16\,d\,\sigma(d)$.
\end{enumerate}
\end{lemma}

One can prove Lemma~\ref{L5.1} applying  the genus $0$ and $1$ TRR
formulas for the descendent constraint $\tau(F)$. In fact, these
TRR formulas consist of:
\begin{enumerate}
\item[(1)] the relation between the tautological class $\psi$ (see below) and
some boundary strata of the Deliegn-Mumford space $\CM_{g,k}$,
\item[(2)] the relation between $\tau(F)$ constraint and $\psi(F)$ constraint.
\end{enumerate}
The relation (1) is also called TRR formula and the relation (2)
follows from {\em relations between generalized correlators}
 (Theorem 1.2 of \cite{km1}). In our case, the computation using
TRR formulas for $\tau(F)$ is quite complicated,
so we will separate the computation into two steps:
we first
use (1) to compute relevant GW invariants of $E(0)$ with
$\psi(F)$ constraints and then apply (2) to those invariants with
$\psi(F)$ constraints to compute invariants with $\tau(F)$ constraints
shown in Lemma~\ref{L5.1}.

In the next section, we  apply the Symplectic sum formula of
\cite{ip3} to the invariants in Lemma~\ref{L5.1} to  compute the
(partial) Gromov--Taubes invariants that appeared in the proof of
Proposition~\ref{P:sum-formula}. After various preliminary lemmas,
we give the proof of Lemma~\ref{L5.1} at the end of this section.

Let $\psi_{i}$ be the  first Chern class of the relative line
bundle $L_{i}$ over $\CM_{g,k}$ whose geometric fiber at the point
$(C;x_{1},\cdots,x_{k})$ is $T^{*}_{x_{i}}$. When $g=0,1$, there
are relations between the class $\psi$ and some boundary strata of
$\CM_{g,k}$ (cf section 4 of \cite{g}). Combining with the
composition law of GW invariants \cite{rt}, those relations give
the TRR formulas for GW invariants of $E(0)$ with $\psi(F)$
constraint: Let $\{H_{\a}\}$ and $\{H^{\a}\}$ be bases of $E(0)$
dual by the intersection form. For
$\b=\b_{1}\otimes\cdots\otimes\b_{n}$ in $[H_{*}(E(0);{\Bbb
Z})]^{\otimes n}$ and an unordered partition of
$\pi=(\pi_{1},\pi_{2})$ of $\{1,\cdots,n\}$ with
$\pi_{1}\ne\emptyset$, we set
$\b_{\pi_{i}}=\b_{l_{1}}\otimes\cdots\otimes \b_{l_{k}}$ where
$\pi_{i}=\{l_{1},\cdots,l_{k}\}$ and $l_{1}<\cdots<l_{k}$. We then
have
\begin{align}
 &\Phi_{A,0}\big(\psi(F),\b_{1},\cdots,\b_{n+2}\big)\notag \\
 &=\,\pm\,
 \sum_{\a}\sum\,
 \Phi_{A_{1},0}\big(F,\b_{\pi_{1}},H_{\a}\big)\,
 \Phi_{A_{2},0}\big(H^{\a},\b_{\pi_{2}},\b_{n+1},\b_{n+2}\big)
 \label{g=0TRR} \\
 &\Phi_{A,1}\big(\psi(F),\b_{1},\cdots,\b_{n}\big)\notag \\
 &=\,\sum_{\a}\tfrac{1}{24}\,
  \Phi_{A,0}\big(F,\b_{1},\cdots,\b_{n},H_{\a},H^{\a}\big) \notag\\
  &\pm\,\sum_{\a} \sum\,
  \Phi_{A_{1},0}\big(F,\b_{\pi_{1}},H_{\a}\big)\,
  \Phi_{A_{2},1}\big(H^{\a},\b_{\pi_{2}}\big)
  \label{g=1TRR}
\end{align}
where the sum is over $A=A_{1}+A_{2}$ and partitions $\pi$ as
above, and the sign depends on the permutation $(\pi_{1},\pi_{2})$
and the degree of $\b_{i}$. In particular, if $\mbox{deg}(\b_{i})$
are all even for $1\leq i\leq n$ the sign is positive.

From now on, we always denote  the fundamental class of $E(0)$ by
1.

\begin{lemma}\label{L5.2}
Let $\Phi$ denote the GW invariants of $E(0)$. Then,
\begin{enumerate}
\item[\rm(a)] $\Phi_{S+dF,1}\big(\psi(F),pt\big)\, =\,2\,\sigma(d)$,
\item[\rm(b)] $\Phi_{S,0} \big(\psi(F),S,F^2\big)\,=\,1$,
\item[\rm(c)] $\Phi_{S,0}\big(\psi(F),pt,F,1\big) =1$,
\item[\rm(d)] $\Phi_{S+dF,1}\big(\psi(F)^2,S\big) \,=\,4\,\sigma(d)$,
\item[\rm(e)] $\Phi_{2S,0} \big(\psi(F)^2,pt,F^2\big) \,=\,2$,
\item[\rm(f)] $\Phi_{dF,1} \big(\psi(F),S,1\big) \,=\,0$.
\end{enumerate}
%
\end{lemma}

\begin{proof}  (a)\qua It follows from  the genus 1 TRR formula
(\ref{g=0TRR}) and Remark~\ref{rat-cur-E(0)} that
\begin{align}\label{L5.2-a1}
 \Phi_{S+dF,1}\big(\psi(F),pt\big)\,&=\,
 \sum_{\a}\tfrac{1}{24}\,\delta_{d,0}\,\Phi_{S,0}\big(F,pt,H_{\a},H^{\a}\big)
 \,+\, \Phi_{S,0}\big(F,pt)\,\Phi_{dF,1}\big(1,pt\big)
 \notag\\
 &+\,
 \Phi_{S,0}\big(F^2,pt\big)\,\Phi_{dF,1}\big(S\big).
\end{align}
The first term in the right hand side vanishes by
Remark~\ref{rat-cur-E(0)}. The second term also vanishes since
$\Phi_{dF,1}\big(pt,\cdot\, \big)=0$. The last term equals
$2\,\sigma(d)$ by Lemma~\ref{well-known}\,a,f. Thus, (a) follows
from (\ref{L5.2-a1}).

\medskip
(b)\qua We have
\begin{equation*}
 \Phi_{S,0}\big(\psi(F),S,F^2\big)\,=\,
 \Phi_{0,0}\big(F,S,1\big)\,
 \Phi_{S,0}\big(pt,F^2\big)\,=\,1
\end{equation*}
where the first equality follows from  the genus $g=0$ TRR formula
(\ref{g=0TRR}) and the second equality follows from
Lemma~\ref{GW=0} and Lemma~\ref{well-known}\,a.

\medskip
(c)\qua We have
\begin{equation*}
 \Phi_{S,0}\big(\psi(F),pt,F,1\big)\,=\,
 \Phi_{S,0}\big(F^2,pt\big) \,
 \Phi_{0,0}\big(S,F,1\big)\,=\,1
\end{equation*}
where the first equality follows from  the genus $g=0$ TRR formula
(\ref{g=0TRR}) and the second equality follows from
Lemma~\ref{GW=0} and Lemma~\ref{well-known}\,a.

\medskip
(d)\qua  It follows from the genus $g=1$ TRR formula (\ref{g=1TRR}) that
\begin{align*}  \label{very-boring-1}
 \Phi_{S+dF,1}\big(\psi(F)^2,S\big)\ &=\ \sum_{\alpha}\tfrac{1}{24} \,
 \Phi_{S+dF,0}\big(F,\psi(F),S,H_{\alpha},H^{\alpha}\big)
 \nonumber \\
 &+\ \sum_{\a}\sum\,
  \Phi_{A_{1},0}\big(F,S,H_{\alpha}\big)\,
  \Phi_{A_{2},1}\big(H^{\alpha},\psi(F)\big)
  \nonumber \\
 &+\ \sum_{\alpha}\sum\,
  \Phi_{A_{1},0}\big(F,\psi(F),H_{\alpha}\big)\,
  \Phi_{A_{2},1}\big(H^{\alpha},S\big)
  \nonumber \\
 &+\,\sum_{\a}\sum\,
  \Phi_{A_{1},0}\big(F,\psi(F),S,H_{\alpha}\big)\,
  \Phi_{A_{2},1}\big(H^{\alpha}\big)
\end{align*}
where the sum is over all decompositions $A_{1}+A_{2}=S+dF$. The
first term in the right hand side vanishes by Remark~\ref
{rat-cur-E(0)}. The second term becomes
\begin{equation*}
 \Phi_{0,0}\big(F,S,1\big)\, \Phi_{S+dF,1}\big(pt,\psi(F)\big).
\end{equation*}
This equals $2\,\sigma(d)$ by Lemma~\ref{GW=0}  and (a).
The third term vanishes since $\CM_{0,3}=\{pt\}$ and the last term
becomes
\begin{equation*}
 \Phi_{S,0}\big(F^2,\psi(F),S\big)\,
 \Phi_{dF,1}\big(S\big).
\end{equation*}
This equals $2\,\sigma(d)$ by  (b) and
Lemma~\ref{well-known}\,f. Thus, we have (d).

\medskip

(e)\qua The genus 0 TRR formula (\ref{g=0TRR}) and Lemma~\ref{GW=0}  give
\begin{align*}
 \Phi_{2S,0}\big(\psi(F)^2,pt,F^2\big)\
 &=\ \sum_{\alpha}
        \Phi_{S,0}\big(F,\psi(F),pt,H_{a}\big)\,
        \Phi_{S,0} \big(H^{a},F^2\big)\notag \\
 &+\   \sum_{\a}  \Phi_{S,0}\big(F,pt,H_{a}\big)\,
        \Phi_{S,0}\big( H^{a},\psi(F),F^2\big)  \notag \\
 &+\ \sum_{\alpha}\,
       \Phi_{S,0}\big(F,\psi(F),H_{a}\big)\,
       \Phi_{S,0}\big( H^{a},pt,F^2\big).
\end{align*}
The first term in the right hand side becomes
\begin{equation*}
 \Phi_{S,0}\big(F,\psi(F),pt,1\big)\,
 \Phi_{S,0} \big(pt,F^2\big).
\end{equation*}
This equals 1 by (c) and Lemma~\ref{well-known}\,a. The second
term becomes
\begin{equation*}
 \Phi_{S,0}\big(F^2,pt\big)\,
 \Phi_{S,0}\big(S,\psi(F),F^2\big).
\end{equation*}
This equals 1 by Lemma~\ref{well-known}\,a and (b).
Since $\CM_{0,3}=\{pt\}$, the last
term vanishes.  Thus, we have (e).

\medskip

(f)\qua The genus $g=1$ TRR formula (\ref{g=1TRR}) and
Remark~\ref{rat-cur-E(0)} give
\begin{align*}
 \Phi_{dF,1}\big(\psi(F),S,1\big)\
 &=\ \sum_{\alpha}\, \tfrac{1}{24}\,\delta_{d0}\,
 \Phi_{0,0}\big( F,S,1,H_{\alpha},H^{\alpha}\big)\\
 &+\
 \sum_{\a}
 \Phi_{0,0}\big(F,S,H_{a}\big)\,\Phi_{dF,1}\big(H^{\alpha},1\big) \\
 &+\  \sum_{\alpha}\,
 \Phi_{0,0}\big(F,1,H_{\alpha}\big)\, \Phi_{dF,1} \big(H^{\alpha},S\big)\\
 &+\
 \sum_{\a}
 \Phi_{0,0}\big(F,1,S,H_{a}\big)\, \Phi_{dF,1} \big(H^{\alpha}\big).
\end{align*}
The first term in the right hand side vanishes by
Lemma~\ref{GW=0}. The second term becomes
\begin{equation*}
 \Phi_{0,0}\big(S,F,1\big)\, \Phi_{dF,1}\big(pt,1\big).
\end{equation*}
This vanishes by the fact $\Phi_{dF,1}\big(pt,\cdot\,\big)=0$. The
third term becomes
\begin{equation*}
 \Phi_{0,0}\big(F,1,S\big)\, \Phi_{dF,1}\big(F,S\big).
\end{equation*}
This also vanishes by the fact
$\Phi_{dF,1}\big(pt,\cdot\,\big)=0$. The last term vanishes as
well by  Lemma~\ref{GW=0}. Thus, we have (f). \end{proof}

Observe that by Lemma~\ref{GW=0}, Remark~\ref{rat-cur-E(0)} and
the dimension count the invariant $\Phi_{A,0}(F,B)$ vanishes
unless $A=S$ and $B=pt$ and that by Lemma~\ref{well-known}\,a and
the Divisor Axiom we have $\Phi_{S,0}(F,pt)=1$. This observation
together with Theorem 1.2 of \cite{km1} gives
\begin{equation}\label{ReCor}
 \Phi_{A,g}\big(\tau(F),\cdot\,\big)\,=\,
 \Phi_{A,g}\big(\psi(F),\cdot\,\big)\,+\,
 \Phi_{A-S,g}\big(1,\cdot\,\big).
\end{equation}

\begin{lemma}\label{L5.3}
Let $\Phi$ denote the GW invariants of $E(0)$. Then,
\begin{enumerate}
\item[\rm(a)] $\Phi_{S,0}\big(\tau(F),F,S\big)\,=\,1$,
\item[\rm(b)] $\Phi_{S,0} \big(\tau(F),pt,F,1\big) =1$.
\end{enumerate}
\end{lemma}

\begin{proof} (a)\qua It follows from (\ref{ReCor}), the fact $\CM_{0,3}=\{pt\}$
and Lemma~\ref{GW=0} that
\begin{equation*}
 \Phi_{S,0}\big(\tau(F),F,S\big)\,=\,
 \Phi_{S,0}\big(\psi(F),F,S\big) \,+\,
 \Phi_{0,0}\big(1,F,S\big)\,=\,1.
\end{equation*}

(b)\qua We have
\begin{equation*}
 \Phi_{S,0}\big(\tau(F),pt,F,1\big)\,=\,
 \Phi_{S,0}\big(\psi(F),pt,F,1\big)\,+\,
 \Phi_{0,0}\big(1,pt,F,1\big)\,=1.
\end{equation*}
where the first equality follows from (\ref{ReCor}) and
the second equality follows from Lemma~\ref{L5.2}\,c
and Lemma~\ref{GW=0}.
\end{proof}

Note that for $B\in H_{2}(E(0);{\Bbb Z})$ the dot product $B\cdot
F\in H_{0}(E(0);\Z)$ corresponds under Poincar\'{e} duality to the
cup product in cohomology. The generalized Divisor Axiom (Lemma
1.4 of \cite{km1}) thus yields
\begin{equation}\label{DA-1}
 \Phi_{A,g}\big(\tau(F),B,\cdot \,\big)\,=\,
 (B\cdot A)\,\Phi_{A,g}\big(\tau(F),\cdot\,\big)
 \,+\,(B\cdot F)\,\Phi_{A,g}\big(pt,\cdot\,\big).
\end{equation}
Combining this relation with (\ref{ReCor}) then gives
\begin{align}\label{DA-2}
 &\Phi_{A,g}\big(\psi(F),B,\cdot\,\big)
 \notag\\
 &=\
 (B\cdot A)\,\Phi_{A,g}\big(\psi(F),\cdot\,\big)
 \,+\,(B\cdot F)\,\Phi_{A,g}\big(pt,\cdot\,\big)
 \,+\,
 (B\cdot S)\,\Phi_{A-S,g}\big(1,\cdot\,\big).
\end{align}

\begin{lemma}\label{L5.4}
Let $\Phi$ denote the GW invariants of $E(0)$. Then,
\begin{enumerate}
\item[\rm(a)] $\Phi_{S,0} \big(\psi(F),1,\g_{1},\g_{2}\big)=1$,
\item[\rm(b)] $\Phi_{2S+dF,1} \big(\psi(F)^3,pt\big)
             \,=\,12\,\sigma(d)$,
\item[\rm(c)] $\Phi_{S+dF,1} \big(\psi(F)^2,1,pt\big)
            \,=\,4\,\sigma(d)$,
\item[\rm(d)] $\Phi_{2S+dF,1}\big(\psi(F),\tau(F),pt^{2}\big) \,=\,
            14\,d\,\sigma(d)$,
\item[\rm(e)] $\Phi_{S+dF,1} \big(\psi(F),1,pt^{2}\big) \,=\,2\,d\,\sigma(d)$.
\end{enumerate}
\end{lemma}

\begin{proof}  (a)\qua It follows from the genus 0 TRR formula (\ref{g=0TRR}),
Lemma~\ref{GW=0} and Lemma~\ref{well-known}\,b that
\begin{equation*}
 2\,\Phi_{S,0} \big(\psi(F),1,\g_{1},\g_{2}\big)\,=\,
 2\,\Phi_{0,0}\big(F,1,S\big)\,
 \Phi_{S,0}\big(F,\g_{1},\g_{2}\big)
 \,=\,2.
\end{equation*}

(b)\qua  The genus $g=1$ TRR formula (\ref{g=1TRR}) and
Remark~\ref{rat-cur-E(0)} give
\begin{align*}
  &\Phi_{2S+dF,1}\big(\psi(F)^3,pt\big)
  \notag \\
  &=\ \sum_{\alpha}\, \tfrac{1}{24}\,\delta_{d0}\,
   \Phi_{2S,0}\big(F,\psi(F)^2,pt,H_{\a},H^{ \a}\big)
   \nonumber \\
  &+\  \sum_{\alpha}\sum\,
  \Phi_{A_{1},0}\big(F,pt,H_{\alpha}\big)\,
  \Phi_{A_{2},1}\big(H^{\alpha},\psi(F)^2\big)
  \nonumber \\
  &+\ \sum_{\alpha}\sum\, 2\,
  \Phi_{A_{1},0}\big(F,\psi(F),H_{\alpha}\big) \,
  \Phi_{A_{2},1}\big(H^{\alpha},\psi(F),pt\big)
  \nonumber \\
  &+\ \sum_{\alpha}\sum\, 2\,
  \Phi_{A_{1},0}\big(F,\psi(F),pt,H_{\alpha} \big)\,
  \Phi_{A_{2},1}\big(H^{\alpha},\psi(F)\big)
  \nonumber \\
  &+\  \sum_{\alpha}\sum\, \,
  \Phi_{A_{1},0}\big(F,\psi(F)^2,H_{ \alpha}\big)\,
  \Phi_{A_{2},1}\big(H^{\alpha},pt\big)
  \nonumber \\
  &+\ \sum_{\alpha}\sum\,
  \Phi_{A_{1},0}\big(F,\psi(F)^2,pt,H_{ \alpha}\big)\,
  \Phi_{A_{2},1}\big(H^{\alpha}\big)  \label{aic-2}
\end{align*}
where the sum is over all decompositions $A_{1}+A_{2}=2S+dF$. The
first term in the right hand side vanishes by
Remark~\ref{rat-cur-E(0)}. The second term becomes
\begin{equation*}
 \Phi_{S,0}\big(F^2,pt\big)\,\Phi_{S+dF,1}\big(S,\psi(F)^2\big).
\end{equation*}
This equals $4\,\sigma(d)$ by Lemma~\ref{well-known}\,a and
Lemma~\ref{L5.2}\,d. The third term vanishes since
$\CM_{0,3}=\{pt\}$. The fourth term becomes
\begin{equation*}
  2\,\Phi_{S,0}\big(F,\psi(F),pt,1\big) \,
  \Phi_{S+dF,1}\big(pt,\psi(F)\big).
\end{equation*}
This equals $4\,\sigma(d)$ by Lemma~\ref{L5.2}\,c,a.
Since $\mbox{dim}_{\Bbb C}\,\CM_{0,4}=1$,
the fifth term vanishes.  The last
term becomes
\begin{equation*}
 \Phi_{2S,0}\big(F^2,\psi(F)^2,pt\big)\,
 \Phi_{dF,1}\big(S\big).
\end{equation*}
by dimension count and Remark~\ref{rat-cur-E(0)}. This equals
$4\,\sigma(d)$ by Lemma~\ref{L5.2}\,e and
Lemma~\ref{well-known}\,f. Thus, we have (b).

\medskip
(c)\qua The genus $g=1$ TRR formula (\ref{g=1TRR}) gives
\begin{align}  \label{L5.4-c1}
 &\Phi_{S+dF,1}\big(\psi(F)^2,1,pt\big) \notag \\
 &=\
 \sum_{\a}\tfrac{1}{24}\,
  \Phi_{S+dF,0}\big(F,\psi(F),1,pt,H_{\a},H^{\a}\big)
 \notag \\
 &+\  \sum_{\a}\sum\,
  \Phi_{A_{1},0}\big(F,1,H_{\alpha}\big)\,
  \Phi_{A_{2},1}\big(H^{\alpha},\psi(F),pt\big)
 \notag \\
 &+\ \sum_{\a}\sum\,
  \Phi_{A_{1},0}\big(F,pt,H_{\alpha}\big)\,
  \Phi_{A_{2},1}\big(H^{\alpha},\psi(F),1,\big)
 \notag \\
 &+\ \sum_{\alpha}\sum\,
  \Phi_{A_{1},0}\big(F,\psi(F),H_{\alpha}\big)\,
  \Phi_{A_{2},1}\big(H^{\alpha},1,pt\big)
 \notag \\
 &+\ \sum_{\a}\sum\,
  \Phi_{A_{1},0}\big(F,1,pt,H_{\a}\big)\,
  \Phi_{A_{2},1}\big(H^{\a},\psi(F)\big)
 \notag \\
 &+\ \sum_{\a}\sum\,
  \Phi_{A_{1},0}\big(F,\psi(F),1,H_{\a}\big)\,
  \Phi_{A_{2},1}\big(H^{\a},pt\big)
 \notag \\
 &+\ \sum_{\a}\sum\,
  \Phi_{A_{1},0}\big(F,\psi(F),pt,H_{\alpha}\big)\,
  \Phi_{A_{2},1}\big(H^{\alpha},1\big)
 \notag \\
 &+\ \sum_{\a}\sum\,
  \Phi_{A_{1},0}\big(F,\psi(F),1,pt,H_{\a}\big)\,
  \Phi_{A_{2},1}\big(H^{\a}\big)
\end{align}
where the sum is over all decompositions $A_{1}+A_{2}=2S+dF$. The
first term in the right hand side of (\ref{L5.4-c1}) vanishes by
Remark~\ref{rat-cur-E(0)}. The second term becomes
\begin{equation}\label{L5.4-c2}
 \Phi_{0,0}\big(F,1,S\big)\,
 \Phi_{S+dF,1}\big(F,\psi(F),pt\big)
\end{equation}
since $\Phi_{dF,1}\big(pt,\cdot\,\big)=0$.  The first factor of
(\ref{L5.4-c2}) is 1 by Lemma~\ref{GW=0} and the second factor is
$2\,\sigma(d)$ by the Divisor Axiom (\ref{DA-2}) and
Lemma~\ref{L5.2}\,a. Thus the second term in the right hand side
of  (\ref{L5.4-c1}) is $2\,\sigma(d)$. The third term vanishes by
Lemma~\ref{GW=0}  and Lemma~\ref{L5.2}\,f;
\begin{align*}
  &\sum_{\a}\sum\,
  \Phi_{A_{1},0}\big(F,pt,H_{\alpha}\big)\,
  \Phi_{A_{2},1}\big(H^{\alpha},\psi(F),1\big)\\
  &=\
  \Phi_{S,0}\big(F^2,pt\big)\,\Phi_{dF,1}\big(S,\psi(F),1\big)\ =\ 0.
\end{align*}
The fourth term vanishes since $\CM_{0,3}=\{pt\}$. The fifth term
vanishes by routine dimension count and Remark~\ref{rat-cur-E(0)}.
The sixth term and the seventh term vanish by routine dimension
count and the fact $\Phi_{dF,1}\big(pt,\cdot\,\big)=0$. The last
term in the right hand side of (\ref{L5.4-c1}) becomes
\begin{equation*}
 \Phi_{S,0}\big(F^2,\psi(F),1,pt\big)\,
 \Phi_{dF,1}\big(S\big).
\end{equation*}
This equals $2\,\sigma(d)$ by Lemma~\ref{L5.2}\,c, the Divisor
Axiom (\ref{DA-2}) and Lemma~\ref{well-known}\,f. Therefore, we
have (c).

\medskip

(d)\qua The genus $g=1$ TRR formula (\ref{g=1TRR}) gives
\begin{align*}
 &\Phi_{2S+dF,1}\big(\psi(F),\tau(F),pt^{2}\big)\notag \\
 &=\
 \sum_{\alpha}\tfrac{1 }{24}\,
 \Phi_{2S+dF,0}\big(F,\tau(F),pt^{2},H_{\alpha},H^{\alpha}\big)
 \nonumber \\
 &=\ \sum_{\alpha}\sum 2\,
 \Phi_{A_{1},0}\big(F,pt,H_{\alpha}\big)\,
 \Phi_{A_{2},1}\big(H^{\alpha},\tau(F),pt\big)
 \nonumber \\
 &=\ \sum_{\alpha}\sum
 \Phi_{A_{1},0}\big(F,\tau(F),H_{\alpha}\big)\,
 \Phi_{A_{2},1}\big(H^{\alpha},pt^{2}\big)
 \nonumber \\
 &=\ \sum_{\alpha}\sum 2\,
 \Phi_{A_{1},0}\big(F,\tau(F),pt,H_{\alpha}\big)\,
 \Phi_{A_{2},1}\big(H^{\alpha},pt\big)
 \nonumber \\
 &=\ \sum_{\alpha}\sum
 \Phi_{A_{1},0}\big(F,pt^{2},H_{\alpha}\big)\,
 \Phi_{A_{2},1}\big(H^{\alpha},\tau(F)\big)
 \nonumber \\
 &=\ \sum_{\alpha}\sum
 \Phi_{A_{1},0}\big(F^2,\tau(F),pt^{2},H_{\alpha}\big)\,
 \Phi_{A_{2},1}\big(H^{\alpha}\big)
\end{align*}
where the sum is over all decompositions $A_{1}+A_{2}=2S+dF$. The
first term of the right hand side vanishes by
Remark~\ref{rat-cur-E(0)}. The second term becomes
\begin{equation*}
 2\,\Phi_{S,0}\big(F^2,pt\big)\,
 \Phi_{S+dF,1}\big(S,\tau(F),pt\big).
\end{equation*}
This equals $8\,d\,\sigma(d)$ by Lemma~\ref {well-known}\,a,
the Divisor Axiom (\ref{DA-1}) and Lemma~\ref {well-known}\,c,d;
\begin{equation*}
 \Phi_{S+dF,1}\big(S,\tau(F),pt\big)\,=\,
 d\,\Phi_{S+dF,1}\big(\tau(F),pt\big)\,+\,
 \Phi_{S+dF,1}\big(pt^2\big)\,=\,
 4\,d\,\sigma(d).
\end{equation*}
The third term becomes
\begin{equation*}
 \Phi_{S,0}\big(F,\tau(F),S\big)\,
 \Phi_{S+dF,1}\big(F,pt^{2}\big)   \,=\,
 2\,d\,\sigma(d).
\end{equation*}
This equals $2\,d\,\sigma(d)$ by Lemma~\ref{L5.3}\,a and
Lemma~\ref{well-known}\,d. The fourth term becomes
\begin{equation*}
 2\,\Phi_{S,0}\big(F,\tau(F),pt,1\big)\,
 \Phi_{S+dF,1}\big(pt^2\big).
\end{equation*}
This equals $4\,d\,\sigma(d)$ by
 Lemma~\ref{L5.3}\,b and  Lemma~\ref{well-known}\,d.
The last two terms vanish by Remark~\ref{rat-cur-E(0)}. Thus, we
have (d).

\medskip
(e)\qua It follows from (\ref{ReCor}) and
 $\Phi_{dF,1}\big(pt,\cdot\,\big)=0$ that
\begin{equation*}
 \Phi_{S+dF,1}\big(\tau(F),1,pt^{2}\big)\,=\,
 \Phi_{S+dF,1}\big(\psi(F),1,pt^{2}\big).
\end{equation*}
By the genus $g=1$ TRR formula (\ref{g=1TRR}) and
Remark~\ref{rat-cur-E(0)}, this becomes
\begin{align*}
 &\sum_{\a}\tfrac{1}{24}\,
   \Phi_{S+dF,1}\big(F,1,pt^{2},H_{\a},H^{\a}\big)\,+\,
 \sum_{\a}
   \Phi_{0,0}\big(F,1,H_{\a}\big)\,
   \Phi_{S+dF,1}\big(H^{\a},pt^{2}\big) \\
 +\,
 &\sum_{\a}2\,
   \Phi_{S,0}\big(F,pt,H_{\a}\big)\,
   \Phi_{dF,1}\big(H^{\a},1,pt\big)\,+\,
 \sum_{\a}2\,
   \Phi_{S,0}\big(F,1,pt,H_{\a}\big)\,
   \Phi_{dF,1}\big(H^{\alpha},pt\big)   \\
 +\,
 &\sum_{\a}
   \Phi_{S,0}\big(F,pt^{2},H_{\a}\big)\,
   \Phi_{dF,1}\big(H^{\a},1\big) \,+\,
 \sum_{\a}
   \Phi_{S,0}\big(F,1,pt^{2},H_{\a}\big)\,
   \Phi_{dF,1}\big(H^{\a}\big).
\end{align*}
The first term vanishes by Remark~\ref{rat-cur-E(0)}. The second
term becomes $2\,d\,\sigma(d)$ by Lemma~\ref{GW=0}  and
Lemma~\ref{well-known}\,d. The third term and the fourth term
vanishes since $\Phi_{dF,1}\big(pt,\cdot\,\big)=0$. The last two
terms vanish by Remark~\ref{rat-cur-E(0)}. Thus, we have (e).
\end{proof}

Now, we are ready to prove Lemma~\ref{L5.1}.

\begin{proof}[Proof of Lemma~\ref{L5.1}]
(a)\qua The relation
(\ref{ReCor}) gives
\begin{equation*}
\Phi_{2S,0}\big(\tau(F)^2,pt\big)\,
 =\,\Phi_{2S,0}\big( \psi(F),\tau(F),pt\big) \,+\,
       \Phi_{S,0}\big(1,\psi(F),pt\big) \,+\,
       \Phi_{0,0}\big(1^2,pt\big).
\end{equation*}
This equals 1 by the fact $\CM_{0,3}=\{pt\}$ and
Lemma~\ref{GW=0}.

\medskip

(b)\qua By (\ref{ReCor}) we have
\begin{align*}
 &\Phi_{2S,0} \big(\tau(F)^2,\g_{1},\g_{2}\big)\\
 &=\, \Phi_{2S,0}\big(\psi(F)^2,\g_{1},\g_{2}\big)\,+\,
 2\,\Phi_{S,0}\big(1,\psi(F),\g_{1},\g_{2}\big) \, +\,
 \Phi_{0,0}\big(1^2,\g_{1},\g_{2}\big).
\end{align*}
Then, (b) follows from the fact $\mbox{dim}_{\Bbb C}\CM_{0,4}=1$,
Lemma~\ref{L5.4}\,a and Lemma~\ref{GW=0}.

\medskip

(c)\qua It follows from (\ref{ReCor}) that
\begin{align*}
 &\Phi_{2S+dF,1}\big(\tau(F)^3,pt\big)\\
 &=\,\Phi_{2S+dF,1}\big(\psi(F)^3,pt\big)\,+\,
 3\,\Phi_{S+dF,1}\big(1,\psi(F)^2,pt\big)
 \, +\,
 3\,\Phi_{dF,1}\big(1^2,\psi(F),pt\big).
\end{align*}
Now, Lemma~\ref{L5.4}\,b,c and the fact
$\Phi_{dF,1}(pt,\cdot\,)=0$ show (c).

\medskip

(d)\qua Using (\ref{ReCor}) yields
\begin{align*}
 &\Phi_{2S+dF,1}\big(\tau(F)^2,pt^{2}\big)\\
 &=\,
  \Phi_{2S+dF,1}\big(\psi(F),\tau(F),pt^{2}\big)\,+\,
  \Phi_{S+dF,1}\big(1,\psi(F),pt^{2}\big)\,+\,
  \Phi_{dF,1}\big(1^2,pt^{2}\big).
\end{align*}
This equals $16\,d\,\sigma(d)$ by Lemma~\ref{L5.4}\,d,e and the
fact  $\Phi_{dF,1}(pt,\cdot\,)=0$.
\end{proof}


\setcounter{equation}{0}
\section{Relative Gromov--Taubes invariants of $E(0)$}
\label{section6}

The aim of this section is to compute the (partial) relative
Gromov--Taubes invariants $G\Phi^{V}$ of $(E(0),V)$ that appeared
in Section \ref{section3}, thereby completing the proof of the sum formulas
(\ref{Sum-1}), (\ref{Sum-2}), (\ref{Sum-1-P}), and
(\ref{Sum-2-P}). Applying the Symplectic Sum Formula of
\cite{ip3}, we will compute those invariants.

The (partial) Gromov--Taubes invariants $G\Phi^{V}$ of $(E(0),V)$
defined in (\ref{GTI}) can be expressed in terms of the relative
invariants $\Phi^{V}$ of $(E(0),V)$. When the multiplicity vector
$s=(2)$, the  invariants $G\Phi^{V}_{s}$ for the classes $2S+dF$,
$d\in {\Bbb Z}$, count $V$--regular maps from a connected domain
and hence
\begin{equation}\label{GT1}
 G\Phi^{V}_{2S+dF,\chi,(2)}\,=\,\Phi^{V}_{2S+dF,g,(2)}\ \ \
 \mbox{where}\ \ \
 g=1-\tfrac{1}{2}\chi.
\end{equation}
When $s=(1,1)$, the invariants $G\Phi^{V}_{s}$ for the classes
$2S+dF$ decompose as a sum of two invariants; one is the invariant
that counts $V$--regular maps from a connected domain, and the
other is the invariant that counts pairs of $V$--regular maps
$(f_{1},f_{2})$, each from a connected domain and having contact
order 1 with $V$. Denote the latter invariants by $T\Phi^{V}$.
Then, we have
\begin{equation}\label{GT2}
 G\Phi^{V}_{2S+dF,\chi,(1,1)}\,=\,
 T\Phi^{V}_{2S+dF,\chi,(1,1)}\,+\,\Phi^{V}_{2S+dF,g,(1,1)}\ \ \
 \mbox{where}\ \ \
 g=1-\tfrac{1}{2}\chi.
\end{equation}
For
 $\beta=\beta_{1}\otimes  \cdots\otimes\beta_{n}\in
 \big[H_{*}\big(E(0)\big)\big]^{\otimes n}$ and
an ordered partition $\pi=(\pi_{1},\pi_{2})$  of $\{1,\cdots,n\}$,
we set
$\beta_{\pi_{i}}=\beta_{l_{1}}\otimes\cdots\otimes\beta_{l_{k}}$
where $\pi_{i}=\{l_{1},\cdots,l_{k}\}$. It then follows  that
\begin{align}
 &T\Phi^{V}_{2S+dF,\chi,(1,1)}
 \big( C_{\g_{i}\cdot\g_{j}};\beta\big)\notag \\
 &=\,\pm\,\sum\,
 \Phi^{V}_{S+d_{1}F,g_{1},(1)} \big(C_{\g_{i}};\beta_{\pi_{1}}\big)
 \, \Phi^{V}_{S+d_{2}F,g_{2},(1)} \big(C_{\g_{j}};\beta_{\pi_{2}}
 \big)  \label{def-T-inv't}
\end{align}
where  the sum is over all $ \chi=(2-2g_{1})+(2-2g_{2})$,
$d=d_{1}+d_{2}$, and partitions $\pi=(\pi_{1},\pi_{2})$ as above,
and the sign depends on the permutation $(\pi_{1},\pi_{2})$ and
the degrees of $\beta_{i}$. In particular, if
$\mbox{deg}(\b_{i})$ are all even the sign is positive.

The following lemma computes the invariants $G\Phi^{V}$  appeared
in the proof of Lemma~\ref{gw-rel(gw)}. The proof of this lemma easily follows from
(\ref{GT2}), (\ref{def-T-inv't}), Remark~\ref{rat-cur-E(0)}, and
Lemma~\ref{well-known}.

\begin{lemma} \label{RI-L3.1}
Let $G\Phi^V$ be the invariants of $(E(0),V)$ defined in
(\ref{GTI}). Then,
\begin{enumerate}
\item[\rm(a)] $G\Phi^{V}_{2S+dF,4,(1,1)}\big(C_{pt^{2}}\big)\,=\,\delta_{d0}$,
\item[\rm(b)] $G\Phi^{V}_{2S+dF,4,(1,1)}\big(C_{pt\cdot F};\g_{1},\g_{2}\big)
           \,=\,\delta_{do}$,
\item[\rm(c)] $G\Phi^{V}_{2S+dF,4,(1,1)}
            \big(C_{(-\g_{1})\cdot \g_{2}};\g_{1},\g_{2}\big)
            \,=\,-\delta_{d0}$,
\item[\rm(d)] $G\Phi^{V}_{2S+dF,2,(2)}\big(C_{pt};\g_{1},\g_{2}\big)
            \,=\,0$,
\item[\rm(e)] $G\Phi^{V}_{2S+dF,2,(1,1)}\big(C_{pt^{2}};\g_{1},\g_{2} \big)
            \,=\,0$,
\item[\rm(f)] $G\Phi^{V}_{2S+dF,4,(1,1)}\big(C_{pt\cdot F};pt\big)
            \,=\,\delta_{d0}$,
\item[\rm(g)] $G\Phi^{V}_{2S+dF,2,(2)}\big(C_{pt};pt\big)
             \,=\,0$,
\item[\rm(h)] $G\Phi^{V}_{2S+dF,2,(1,1)}\big(C_{pt^{2}};pt\big)
             \,=\,2\,d\,\sigma(d)$.
\end{enumerate}
\end{lemma}

The following lemma lists the invariants $G\Phi^{V}$ that entered
in the proof of Proposition~\ref{P:sum-formula}\,a and
is an immediate consequence of (\ref{GT2}), (\ref{def-T-inv't}),
Remark~\ref{rat-cur-E(0)}, Lemma~\ref{well-known}, and
Lemma~\ref{lemma5-1}.

\begin{lemma}
\label{RI-P3.3a} Let $G\Phi^V$ be the invariants of $(E(0),V)$
defined in (\ref{GTI}). Then,
\begin{enumerate}
\item[\rm(a)] $G\Phi^{V}_{2S+dF,4,(1,1)}\big(C_{pt\cdot F};\tau(F)\big)
            \,=\,0$,
\item[\rm(b)] $G\Phi^{V}_{2S+dF,2,(2)}\big(C_{pt};\tau(F)\big)\,=\,
            \frac{1}{2}\,\delta_{d0}$,
\item[\rm(c)] $G\Phi^{V}_{2S+dF,2,(1,1)}\big(C_{pt^{2}};\tau(F)\big)
             \,=\,4\,\sigma(d)$.
\end{enumerate}
\end{lemma}

Lastly, Lemma~\ref{RI-P3.3b} below computes
the invariants $G\Phi^V$ that appeared in the proof of
 Proposition~\ref{P:sum-formula}\,b.
In order to prove this lemma, we will apply the Symplectic Sum
Formula of \cite{ip3} to the GW invariants of $E(0)$ in
shown Lemma~\ref{L5.1} by writing $E(0)$ as a symplectic sum
\begin{equation}\label{SS-E(0)}
   E(0)\,=\,E(0)\,\#_{V}\,E(0)
\end{equation}
and by splitting constraints in various ways. In this case, there
is also no contribution from the neck (\,cf Lemma 16.1 of
\cite{ip3}) and hence we have the following sum formulas: Let
$\{\g_{i}\}$  be a basis of $H_{*}(V;{\Bbb Z})$ and  $\{\g^{i}\}$
be its dual basis with respect to  the intersection form of $V$.
For $\b=\b_{1}\otimes\cdots\otimes\b_{n}$ in $[H_{*}(E(0);{\Bbb
Z})]^{\otimes n}$ and $n=n_{1}+n_{2}$, we set
\begin{equation*}
 \b^\prime=\b_{1}\otimes\cdots\otimes \b_{n_{1}}\ \ \ \
 \mbox{and}\ \ \ \
  \b^{\prime\prime}=\b_{n_{1}+1}\otimes\cdots\otimes \b_{n}.
\end{equation*}
For $\tau(F)^m$ and $m=m_1+m_2$, we set
\begin{equation*}
\tau^\prime\ =\ \tau(F)^{m_{1}}\ \ \ \
\mbox{and}\ \ \ \
\tau^{\prime\prime}\ =\ \tau(F)^{m_2}.
\end{equation*}
Then, for such splitting of constraints $\b=\b^{\prime}\otimes
\b^{\prime\prime}$ and
$\tau(F)^m=\tau^{\prime}\cdot\tau^{\prime\prime}$ the sum formula of
the symplectic sum (\ref{SS-E(0)}) for the class $S+dF$ becomes
\begin{align}\label{S+dF}
  &\Phi_{S+dF,g}\big(\b,\tau(F)^m\big)\notag \\
  &=\ \sum\,\sum_{i}
   \Phi^{V}_{S+d_{1}F,g_{1},(1)}
     \big(\b^{\prime},\tau^{\prime};C_{\g_{i}}\big)\,
   \Phi^{V}_{S+d_{2}F,g_{2},(1)}
     \big(C_{\g^{i}};\b^{\prime\prime},\tau^{\prime\prime}\big)
\end{align}
where the sum is over all $d=d_{1}+d_{2}$ and $g=g_{1}+g_{2}$.
Similarly, the sum formula of the sum (\ref{SS-E(0)})
applied to the class $2S+dF$ with
the splitting of constraints $\b=\b^{\prime}\otimes
\b^{\prime\prime}$ and $\tau(F)^m=\tau^{\prime}\cdot\tau^{\prime\prime}$
gives
\begin{align}\label{2S+dF}
  &\Phi_{2S+dF,g}\big(\b,\tau(F)^m\big)
  \notag \\
  &=\ \sum\,\sum_{i,j} \tfrac{1}{2}\,
  \Phi^{V}_{2S+d_{1}F,g_{1},(1,1)}
    \big(\b^{\prime},\tau^{\prime};C_{\g_{i}\cdot \g_{j}}\big)\,
  T\Phi^{V}_{2S+d_{2}F,\chi_{2},(1,1)}
    \big(C_{\g^{j}\cdot\g^{i}};\b^{\prime\prime},\tau^{\prime\prime}\big)
  \notag \\
  &+\ \sum\,\sum_{i,j} \tfrac{1}{2}\,
  \Phi^{V}_{2S+d_{1}F,g_{1},(1,1)}
    \big(\b^{\prime},\tau^{\prime};C_{\g_{i}\cdot \g_{j}}\big)\,
  \Phi^{V}_{2S+d_{2}F,g_{2},(1,1)}
    \big(C_{\g^{j}\cdot\g^{i}};\b^{\prime\prime},\tau^{\prime\prime}\big)
  \notag \\
  &+\ \sum\,\sum_{i} 2\,
  \Phi^{V}_{2S+d_{1}F,g_{1},(2)}
    \big(\b^{\prime},\tau^{\prime};C_{\g_{i}}\big)\,
  \Phi^{V}_{2S+d_{2}F,g_{2},(2)}
    \big(C_{\g^{i}};\b^{\prime},\tau^{\prime\prime}\big)
  \notag \\
  &+\ \sum\,\sum_{i,j} \tfrac{1}{2}\,
  T\Phi^{V}_{2S+d_{1}F,\chi_{1},(1,1)}
    \big(\b^{\prime},\tau^{\prime};C_{\g_{i}\cdot \g_{j}}\big)\,
  \Phi^{V}_{2S+d_{2}F,g_{2},(1,1)}
    \big(C_{\g^{j}\cdot\g^{i}};\b^{\prime\prime},\tau^{\prime}\big)
\end{align}
where the sum is over all $d=d_{1}+d_{2}$, and
$g=g_{1}+2-\frac{1}{2}\chi_{2}$ for the first term,
$g=g_{1}+g_{2}+1$ for the second term,
$g=g_{1}+g_{2}$ for the third term  and
$g=2-\frac{1}{2}\chi_{1}+g_{2}$ for the last term.

\begin{lemma}
\label{RI-P3.3b} Let $G\Phi^V$ be the invariants of $(E(0),V)$
defined in (\ref{GTI}). Then
\begin{enumerate}
\item[\rm(a)]
   $G\Phi^{V}_{2S+dF,4,(1,1)} \big(C_{F^{2}};\tau(F)^2 \big)\,=\,0$,
\item[\rm(b)]
  $G\Phi^{V}_{2S+dF,2,(2)} \big(C_{F};\tau(F)^2\big)\,=\,0$,
\item[\rm(c)]
   $G\Phi^{V}_{2S+dF,2,(1,1)} \big(C_{pt\cdot F};\tau(F)^2\big)\,=\,\delta_{d0}$,
\item[\rm(d)]
  $G\Phi{V}_{2S+dF,2,(1,1)}
    \big(C_{(-\g_{1})\cdot  \g_{2}};\tau(F)^2\big)\,=\,\delta_{d0}$,
\item[\rm(e)]
  $G\Phi^{V}_{2S+dF,0,(2)} \big(C_{pt};\tau(F)^2\big)\,=\,10\,\sigma(d)$,
\item[\rm(f)]
  $ G\Phi^{V}_{2S+dF,0,(1,1)}
      \big(C_{pt^{2}};\tau(F)^2 \big)\,=\,
      \underset{d_{1}+d_{2}=d}{\sum}
      16\,\sigma(d_{1})\,\sigma(d_{2})\,+\,12\,d\,\sigma(d)$.
\end{enumerate}
\end{lemma}

\begin{proof} (a)\qua Using (\ref{GT2}), (\ref{def-T-inv't}),
Remark~\ref{rat-cur-E(0)}, and Lemma~\ref{lemma5-1}\,b gives
\begin{equation*}
  G\Phi^{V}_{2S+dF,4,(1,1)} \big(C_{F^{2}};\tau(F)^2\big)\,=\,
  2\,\delta_{d0}\,
  \Phi^{V}_{S,0}\big( C_{F};\tau(F)\big)\,
  \Phi^{V}_{S,0}\big( C_{F};\tau(F)\big)\,=\,0.
\end{equation*}

(b)\qua By (\ref{GT1}), Remark~\ref{rat-cur-E(0)} and
Lemma~\ref{lemma5-1}\,b we have
\begin{equation*}
  G\Phi^{V}_{2S+dF,2,(2)} \big(C_{F};\tau(F)^2\big)\,=\,\delta_{d0}\,
  \Phi^{V}_{2S,0,(2)} \big(C_{F};\tau(F)^2\big)\,=\,0.
\end{equation*}

(c)\qua It follows from the sum formula
Remark~\ref{sd} that
\begin{equation*}
  \Phi_{S+dF,1}\big(\tau(F)^2\big)\,=\,
  \sum_{d_{1}+d_{2}=d}
  \Phi_{S+d_{1}F,0,(1)}^{V}\big(C_{pt}\big)\,
  \Phi_{S+d_{2}F,1,(1)}^{V} \big(C_{F};\tau(F)^2\big).
\end{equation*}
By Remark~\ref{rat-cur-E(0)} and Lemma~\ref{well-known}\,a,  the
right hand side of this becomes
\begin{equation*}
    \Phi_{S+dF,1,(1)}^{V} \big(C_{F};\tau(F)^2\big).
\end{equation*}
Similarly, the sum formula (\ref{S+dF}), Remark~\ref{rat-cur-E(0)}
and Lemma~\ref{lemma5-1}\,b yields
\begin{align*}
  \Phi_{S+dF,1}\big(\tau(F)^2\big)\,&= \sum_{d_{1}+d_{2}=d}
  \Phi_{S+d_{1}F,1,(1)}^{V}\big(\tau(F);C_{pt}\big)\,
  \Phi_{S+d_{2}F,0,(1)}^{V}\big(C_{F};\tau(F)\big)  \nonumber \\
  &+\sum_{d_{1}+d_{2}=d}
  \Phi_{S+d_{1}F,0,(1)}^{V}\big(\tau(F);C_{F}\big)\,
  \Phi_{S+d_{2}F,1,(1)}^{V}\big(C_{pt};\tau(F)\big)
  \notag \\
  &=\,0.
\end{align*}
Thus, we have
\begin{equation}\label{c0}
  \Phi_{S+dF,1,(1)}^{V} \big(C_{F};\tau(F)^2\big)\,=\,
   \Phi_{S+dF,1}\big(\tau(F)^2\big)\,=\,0.
\end{equation}
This together with  (\ref{def-T-inv't}) and
Lemma~\ref{lemma5-1}\,b then implies that
\begin{align}\label{c1}
  T\Phi_{2S+dF,2,(1,1)}^{V} \big(C_{F\cdot pt};\tau(F)^2\big)\,
  &=\,
  \Phi_{S+dF,1,(1)}^{V} \big(C_{F};\tau(F)^2\big)\,
  \Phi_{S,0,(1)}^{V}  \big(C_{pt}\big) \notag \\
  &+\,
  2\, \Phi_{S,0,(1)}^{V}   \big(C_{F};\tau(F)\big)\,
  \Phi_{S+dF,1,(1)}^{V} \big(C_{pt};\tau(F)\big) \notag \\
  &=\,0.
\end{align}
On the other hand, applying the sum formula (\ref{2S+dF}) we
obtain
\begin{align*}
  \Phi_{2S,0}\big(pt,\tau(F)^2\big)\,
  &=\,\tfrac{1}{2}\,
  \Phi_{2S,0,(1,1)}^{V}\big(pt;C_{pt^{2}}\big)\,
  T\Phi_{2S,4,(1,1)}^{V}\big(C_{F^{2}};\tau(F)^2\big)   \\
 &+\,2\,
 \Phi_{2S,0,(2)}^{V}\big(pt;C_{pt}\big)\,
 \Phi_{2S,0,(2)}^{V}\big(C_{F};\tau(F)^2\big)   \\
 &+\,
 T\Phi_{2S,4,(1,1)}^{V}\big(pt;C_{pt\cdot F}\big)\,
 \Phi_{2S,0,(1,1)}^{V}\big(C_{F\cdot pt};\tau(F)^2\big).
\end{align*}
By Remark~\ref{rat-cur-E(0)}, (\ref{GT2}) and
Lemma~\ref{RI-L3.1}\,f, the right hand side of this can be
simplified as
\begin{equation*}
  \Phi_{2S,0,(1,1)}^{V}\big(C_{F\cdot pt};\tau(F)^2\big).
\end{equation*}
Consequently, we have
\begin{equation}\label{c2}
  \Phi_{2S,0,(1,1)}^{V}\big(C_{F\cdot pt};\tau(F)^2\big)\,=\,
  \Phi_{2S,0}\big(pt,\tau(F)^2\big)\,=\,1
\end{equation}
where the second equality follows from Lemma~\ref{L5.1}\,a. Now,
(c) follows from (\ref{GT2}), Remark~\ref{rat-cur-E(0)},
(\ref{c1}) and (\ref{c2}):
\begin{align*}
  &G\Phi^{V}_{2S+dF,2,(1,1)} \big(C_{pt\cdot F};\tau(F)^2\big) \\
  &=\
  T\Phi^{V}_{2S+dF,2,(1,1)} \big(C_{pt\cdot F};\tau(F)^2\big)\,+\,
  \delta_{d0}\,\Phi^{V}_{2S,0,(1,1)} \big(C_{pt\cdot F};\tau(F)^2\big)\\
  &=\ \delta_{d0}.
\end{align*}

 (d)\qua We have
\begin{align*}
  &\Phi_{2S,0}\big(\g_{1},\g_{2},\tau(F)^2\big)
  \notag \\
  &=\ \tfrac{1}{2}\,
  \Phi_{2S,0,(1,1)}^{V}\big(\g_{1},\g_{2};C_{pt^{2}}\big)\,
  T\Phi_{2S,4,(1,1)}^{V}\big(C_{F^{2}};\tau(F)^2\big)
  \notag   \\
  &+\ 2\,
  \Phi_{2S,0,(2)}^{V}\big(\g_{1},\g_{2};C_{pt}\big)\,
  \Phi_{2S,0,(2)}^{V}\big(C_{F};\tau(F)^2\big)
  \nonumber \\
  &+\
  T\Phi_{2S,4,(1,1)}^{V}\big(\g_{1},\g_{2};C_{pt\cdot F}\big)\,
  \Phi_{2S,0,(1,1)}^{V}\big(C_{F\cdot pt};\tau(F)^2\big)
  \nonumber \\
  &+\
  T\Phi_{2S,4,(1,1)}^{V}\big(\g_{1},\g_{2};C_{\g_{1}\cdot \g_{2}}\big)\,
  \Phi_{2S,0,(1,1)}^{V}\big(C_{(-\g_{1})\cdot \g_{2}};\tau(F)^2\big)
  \\
  &=\
 \Phi_{2S,0,(1,1)}^{V}\big(C_{F\cdot pt};\tau(F)^2\big) \,+\,
 \Phi_{2S,0,(1,1)}^{V}\big(C_{(-\g_{1})\cdot \g_{2}};\tau(F)^2\big).
\end{align*}
where the first equality follows from the sum formula (\ref{2S+dF})
and the second follows from Remark~\ref{rat-cur-E(0)},
(\ref{GT2}) and Lemma~\ref{RI-L3.1}\,b,c.
This together with  Lemma~\ref{L5.1}\,b and (\ref{c2}) shows
$$ \Phi_{2S,0,(1,1)}^{V}
   \big(C_{(-\g_{1})\cdot \g_{2}};\tau(F)^2\big)\,=\,1.$$
On the other hand, (\ref{def-T-inv't}) and
dimension count give
$$ T\Phi^{V}_{2S+dF,2,(1,1)}
     \big(C_{(-\g_{1})\cdot  \g_{2}};\tau(F)^2\big)\,=\,0.$$
Consequently, we have
\begin{align*}
  &G\Phi^{V}_{2S+dF,2,(1,1)}
     \big(C_{(-\g_{1})\cdot  \g_{2}};\tau(F)^2\big) \\
  &=\
  T\Phi^{V}_{2S+dF,2,(1,1)}
     \big(C_{(-\g_{1})\cdot  \g_{2}};\tau(F)^2\big)\,+\,
  \delta_{d0}\,
   \Phi_{2S,0,(1,1)}^{V}\big(C_{(-\g_{1})\cdot \g_{2}};\tau(F)^2\big)\\
   &=\ \delta_{d0}.
\end{align*}

(e)\qua Using the sum formula (\ref{2S+dF}), we have
\begin{align}  \label{s-g-1}
  &\Phi_{2S+dF,1} \big(pt,\tau(F)^3\big)
   \notag \\
  &=\ \sum_{d_{1}+d_{2}=d}\tfrac{1}{2}\,
  \Phi^{V}_{2S+d_{1}F,1,(1,1)}\big( pt,\tau(F);C_{pt^{2}}\big)\,
  T\Phi^{V}_{2S+d_{2}F,4,(1,1)}\big(C_{F^{2}};\tau(F)^2\big)
  \notag \\
  &+\ \sum_{d_{1}+d_{2}=d} 2\,
  \Phi^{V}_{2S+d_{1}F,1,(2)}\big(pt,\tau(F);C_{pt} \big)\,
  \Phi^{V}_{2S+d_{2}F,0,(2)}\big(C_{F};\tau(F)^2\big)
  \notag \\
  &+\ \sum_{d_{1}+d_{2}=d}
  \Phi^{V}_{2S+d_{1}F,0,(1,1)} \big(pt,\tau(F);C_{pt\cdot F}\big)\,
  T\Phi^{V}_{2S+d_{2}F,2,(1,1)} \big(C_{pt\cdot F};\tau(F)^2\big)
  \notag \\
  &+\ \sum_{d_{1}+d_{2}=d}
  \Phi^{V}_{2S+d_{1}F,0,(1,1)} \big(pt,\tau(F);C_{pt\cdot F}\big)\,
  \Phi^{V}_{2S+d_{2}F,0,(1,1)} \big(C_{pt\cdot F};\tau(F)^2\big)
  \notag \\
  &+\ \sum_{d_{1}+d_{2}=d} 2\,
  \Phi^{V}_{2S+d_{1}F,0,(2)} \big(pt,\tau(F);C_{F} \big)\,
  \Phi^{V}_{2S+d_{2}F,1,(2)} \big(C_{pt};\tau(F)^2\big)
  \notag \\
  &+\ \sum_{d_{1}+d_{2}=d}
  T\Phi^{V}_{2S+d_{1}F,2,(1,1)} \big(pt,\tau(F);C_{pt\cdot F}\big)\,
  \Phi^{V}_{2S+d_{2}F,0,(1,1)} \big(C_{pt\cdot F};\tau(F)^2\big)
  \notag \\
  &+\ \sum_{d_{1}+d_{2}=d}\tfrac{1}{2}\,
  T\Phi^{V}_{2S+d_{1}F,4,(1,1)} \big(pt,\tau(F);C_{F^{2}}\big)\,
  \Phi^{V}_{2S+d_{2}F,1,(1,1)} \big(C_{pt^{2}};\tau(F)^2\big).
\end{align}
The first term in the right hand side vanishes by
(\ref{def-T-inv't}) and Lemma~\ref{lemma5-1}\,b. The second term
vanishes by Lemma~\ref{lemma5-1}\,b. The third and the fourth
terms vanish by Remark~\ref{rat-cur-E(0)}. The fifth term becomes
\begin{equation*}
 2\,\Phi^{V}_{2S+dF,1,(2)} \big(C_{pt};\tau(F)^2\big)
\end{equation*}
by Remark~\ref{rat-cur-E(0)} and Lemma~\ref{lemma5-1}\,a. The
sixth term becomes
\begin{equation*}
  T\Phi^{V}_{2S+dF,2,(1,1)} \big(pt,\tau(F);C_{pt\cdot F}\big)
\end{equation*}
by Remark~\ref{rat-cur-E(0)} and (\ref{c2}). The last term
vanishes by (\ref{def-T-inv't}) and Lemma~\ref{lemma5-1}\,b. Thus,
together with Lemma~\ref{L5.1}\,c, we  have
\begin{equation}\label{e1}
 24\,\sigma(d)\,=\,
 2\,\Phi^{V}_{2S+dF,1,(2)} \big(C_{pt};\tau(F)^2\big)\,+\,
 T\Phi^{V}_{2S+dF,2,(1,1)} \big(pt,\tau(F);C_{pt\cdot F}\big).
\end{equation}
On the other hand, it follows from (\ref{def-T-inv't}),
Remark~\ref{rat-cur-E(0)}, Lemma~\ref{well-known}\,a,c and
Lemma~\ref{lemma5-1}\,b that
\begin{align}\label{e2}
  T\Phi^{V}_{2S+dF,2,(1,1)} \big(pt,\tau(F);C_{pt\cdot F}\big)
  &=\
  \Phi^{V}_{S,0,(1)}\big(C_{pt}\big)\,
  \Phi^{V}_{S+dF,1,(1)}\big(pt,\tau(F);C_{F}\big) \notag \\
  &+
  \Phi^{V}_{S+dF,1,(1)}\big(\tau(F);C_{pt}\big)
  \Phi^{V}_{S,0,(1)}\big(pt;C_{F}\big)
  \notag\\
  &+\
  \Phi^{V}_{S+dF,1,(1)}\big(pt;C_{pt}\big)\,
  \Phi^{V}_{S,0,(1)}\big(\tau(F);C_{F}\big)
  \notag\\
  &=\  4\,\sigma(d).
\end{align}
Then, (e) follows from (\ref{GT1}), (\ref{e1}) and (\ref{e2});
\begin{equation*}
  G\Phi^{V}_{2S+dF,0,(2)} \big(C_{pt};\tau(F)^2\big)\,=\,
  \Phi^{V}_{2S+dF,1,(2)} \big(C_{pt};\tau(F)^2\big)\,=\,
  10\,\sigma(d).
\end{equation*}

(f)\qua The sum formula (\ref{2S+dF}) gives
\begin{align}  \label{S-f}
 &\Phi_{2S+dF,1} \big(pt^{2},\tau(F)^2\big)
 \notag \\
 &=\ \sum_{d_{1}+d_{2}=d}\tfrac{1}{2}\,
 \Phi_{2S+d_{1}F,1,(1,1)}^{V}\big(pt^{2};C_{pt^{2}}\big)\,
 T\Phi_{2S+d_{2},4,(1,1)}^{V}\big(C_{F^{2}};\tau(F)^2\big)
 \notag \\
 &+\ \sum_{d_{1}+d_{2}=d}2\,
 \Phi_{2S+d_{1}F,1,(2)}^{V}\big(pt^{2};C_{pt}\big)\,
 \Phi_{2S+d_{2}F,0,(2)}^{V}\big(C_{F};\tau(F)^2\big)
 \notag\\
 &+\ \sum_{d_{1}+d_{2}=d}
 \Phi_{2S+d_{1}F,0,(1,1)}^{V}\big(pt^{2};C_{pt\cdot F}\big)\,
 T\Phi_{2S+d_{2}F,2,(1,1)}^{V}\big(C_{pt\cdot F};\tau(F)^2\big)
 \notag \\
 &+\ \sum_{d_{1}+d_{2}=d}
 \Phi_{2S+d_{1}F,0,(1,1)}^{V}\big(pt^{2};C_{pt\cdot F}\big)\,
 \Phi_{2S+d_{2}F,0,(1,1)}^{V}\big(C_{pt\cdot F};\tau(F)^2\big)
 \notag \\
 &+\ \sum_{d_{1}+d_{2}=d}2\,
 \Phi_{2S+d_{1}F,0,(2)}^{V}\big(pt^{2};C_{F}\big)\,
 \Phi_{2S+d_{2}F,1,(2)}^{V}\big(C_{pt};\tau(F)^2\big)
 \notag\\
 &+\ \sum_{d_{1}+d_{2}=d}
 T\Phi_{2S+d_{1}F,2,(1,1)}^{V}\big(pt^{2};C_{pt\cdot F}\big)\,
 \Phi_{2S+d_{2}F,0,(1,1)}^{V}\big(C_{pt\cdot F};\tau(F)^2\big)
 \notag \\
 &+\ \sum_{d_{1}+d_{2}=d}\tfrac{1}{2}\,
 T\Phi_{2S+d_{1}F,4,(1,1)}^{V}\big(pt^{2};C_{F^{2}}\big)\,
 \Phi_{2S+d_{2}F,1,(1,1)}^{V}\big(C_{pt^{2}};\tau(F)^2\big).
\end{align}
The first two terms in the right hand side vanish by
Lemma~\ref{lemma5-1}\,b, while the next three terms vanish by
Remark~\ref {rat-cur-E(0)}. Thus, by Lemma~\ref{L5.1}\,d and
Remark~\ref{rat-cur-E(0)}, the equation (\ref{S-f}) becomes
\begin{align}\label{S-f2}
 16\,d\,\sigma(d)\,&=\,
 T\Phi_{2S+dF,2,(1,1)}^{V}\big(pt^{2};C_{pt\cdot F}\big)\,
 \Phi_{2S,0,(1,1)}^{V}\big(C_{pt\cdot F};\tau(F)^2\big)
 \notag \\
 &+\,\tfrac{1}{2}\,
 T\Phi_{2S,4,(1,1)}^{V}\big(pt^{2};C_{F^{2}}\big)\,
 \Phi_{2S+dF,1,(1,1)}^{V}\big(C_{pt^{2}};\tau(F)^2\big).
\end{align}
By (\ref{def-T-inv't}) and Lemma~\ref{well-known}\,a,d, the first
$T\Phi^{V}$ invariant in the right hand side becomes
\begin{align}\label{TP-f1}
 &T\Phi_{2S+dF,2,(1,1)}^{V}\big(pt^{2};C_{pt\cdot F}\big)
 \notag \\
 &=\
 \Phi_{S,0,(1)}^{V}\big(C_{pt}\big)\,
 \Phi_{S+dF,1,(1)}^{V}\big(pt^{2};C_{F}\big)\,+\,
 2\,\Phi_{S+dF,1,(1)}^{V}\big(pt;C_{pt}\big)\,
 \Phi_{S,0,(1)}^{V}\big(pt;C_{F}\big)
 \notag \\
 &=\ 4\,d\,\sigma(d).
\end{align}
By (\ref{def-T-inv't}) and Lemma~\ref{well-known}\,a, the second
$T\Phi^{V}$ invariant becomes
\begin{equation}  \label{TP-f2}
 T\Phi_{2S,4,(1,1)}^{V}\big(pt^{2};C_{F^{2}}\big)\,=\,
 2\,\Phi_{S,0,(1)}^{V}\big(pt;C_{F}\big)\,
 \Phi_{S,0,(1)}^{V}\big(pt;C_{F}\big)\,=\,2.
\end{equation}
Then, by (\ref{S-f2}), (\ref{TP-f1}), (\ref{c2}) and (\ref{TP-f2})
we have
\begin{equation}\label{f1}
 \Phi_{2S+dF,1,(1,1)}^{V}\big(C_{pt^{2}};\tau(F)^2\big)
 \,=\,12\,d\,\sigma(d).
\end{equation}

On the other hand, the sum formula (\ref{S+dF}) and
Remark~\ref{rat-cur-E(0)} give
\begin{align}  \label{L-TP-1}
 &\Phi_{S+dF,2}\big(\tau(F)^2\big)
 \notag \\
 &=\ \sum_{d_{1}+d_{2}=d}
 \Phi^{V}_{S+d_{1}F,1,(1)}\big(pt;C_{pt}\big)\,
 \Phi^{V}_{S+d_{2}F,1,(1)}\big(C_{F};\tau(F)^2\big)
 \notag \\
 &+\
 \Phi^{V}_{S,0,(1)}\big(pt;C_{F}\big)\,
 \Phi^{V}_{S+F,2,(1)}\big(C_{pt};\tau(F)^2\big).
\end{align}
The first term of the right hand side vanishes by (\ref{c0}) and
hence by Lemma~\ref{well-known}\,a, the equation (\ref{L-TP-1})
becomes
\begin{equation}\label{L-TP-2}
 \Phi_{S+dF,2}\big(\tau(F)^2\big)\,=\,
 \Phi^{V}_{S+F,2,(1)}\big(C_{pt};\tau(F)^2\big).
\end{equation}
Splitting constraints in a different way,  the sum formula
(\ref{S+dF}) and Remark~\ref{rat-cur-E(0)} give
\begin{align}  \label{L-TP-3}
 &\Phi_{S+dF,2}\big(\tau(F)^2\big)
 \nonumber \\
 &=\
 \Phi^{V}_{S+dF,2,(1)}\big(pt,\tau(F);C_{pt}\big)\,
 \Phi^{V}_{S,0,(1)}\big(C_{F};\tau(F)\big)
 \notag \\
 &+\ \sum_{d_{1}+d_{2}=d}
 \Phi^{V}_{S+d_{1}F,1,(1)}\big(pt,\tau(F);C_{F}\big)\,
 \Phi^{V}_{S+d_{2}F,(1)}\big(C_{pt};\tau(F)\big).
\end{align}
The first term of the right hand side vanishes by
Lemma~\ref{lemma5-1}\,b. Thus, by Lemma~\ref{well-known}\,c the
equation (\ref{L-TP-3}) becomes
\begin{equation}\label{L-TP-4}
 \Phi_{S+dF,2}\big(\tau(F)^2\big)\,=\,
 \sum_{d_{1}+d_{2}=d}4\,\sigma(d_{1})\,\sigma(d_{2}).
\end{equation}
Relating (\ref{L-TP-2}) and (\ref{L-TP-4}), we have
\begin{equation}\label{L-TP-5}
 \Phi^{V}_{S+F,2,(1)}\big(C_{pt};\tau(F)^2\big)\,=\,
 \sum_{d_{1}+d_{2}=d}4\,\sigma(d_{1})\,\sigma(d_{2}).
\end{equation}
Consequently, (\ref{def-T-inv't}), Remark~\ref{rat-cur-E(0)},
(\ref{L-TP-5}) and  Lemma~\ref{well-known}\,a,c show that
\begin{align}\label{f2}
 &T\Phi^{V}_{2S+dF,0,(1,1)}
   \big(C_{pt^{2}};\tau(F)^2\big)\
 \notag \\
 =&\
 2\,\Phi^{V}_{S+dF,2,(1)}\big(C_{pt};\tau(F)^2\big)\,
 \Phi^{V}_{S,0,(1)}\big(C_{pt}\big)
 \notag \\
 +&\ \sum_{d_{1}+d_{2}=d}
 2\,\Phi^{V}_{S+d_{1}F,1,(1)}\big(C_{pt};\tau(F)\big)\,
 \Phi^{V}_{S+d_{2}F,1,(1)}\big(C_{pt};\tau(F)\big)
 \notag \\
 =&\ \sum_{d_{1}+d_{2}=d} 16\,\sigma(d_{1})\,\sigma(d_{2}).
\end{align}
Thus, we have
\begin{align*}
 &G\Phi^{V}_{2S+dF,0,(1,1)}
      \big(C_{pt^{2}};\tau(F)^2 \big)
 \notag \\
 &=\
 T\Phi^{V}_{2S+dF,0,(1,1)}
      \big(C_{pt^{2}};\tau(F)^2 \big)\,+\,
 \Phi^{V}_{2S+dF,1,(1,1)}
      \big(C_{pt^{2}};\tau(F)^2 \big)
 \notag \\
 &=\ \sum_{d_{1}+d_{2}=d} 16\,\sigma(d_{1})\,\sigma(d_{2})\,+\,
 12\,d\,\sigma(d)
\end{align*}
where the first equality follows from (\ref{GT2}) and
the second follows from (\ref{f1}) and (\ref{f2}).
\end{proof}

\end{document}